
\def\LaTeX{%
  \let\Begin\begin
  \let\End\end
  \let\salta\relax
  \let\finqui\relax
  \let\futuro\relax}

\def\UK{\def\our{our}\let\sz s}
\def\USA{\def\our{or}\let\sz z}

\UK



\LaTeX

\USA


\salta

\documentclass[twoside,12pt]{article}
\setlength{\textheight}{24cm}
\setlength{\textwidth}{16cm}
\setlength{\oddsidemargin}{2mm}
\setlength{\evensidemargin}{2mm}
\setlength{\topmargin}{-15mm}
\parskip2mm


\usepackage[usenames,dvipsnames]{color}
\usepackage{amsmath}
\usepackage{amsthm}
\usepackage{amssymb}
\usepackage[mathcal]{euscript}


\usepackage[T1]{fontenc}
\usepackage[latin1]{inputenc}
\usepackage[english]{babel}
\usepackage[babel]{csquotes}

\usepackage{cite}

\usepackage{latexsym}
\usepackage{graphicx}
\usepackage{mathrsfs}
\usepackage{mathrsfs}
\usepackage{hyperref}
\usepackage{pgfplots}


%
%


\definecolor{viola}{rgb}{0.3,0,0.7}
\definecolor{ciclamino}{rgb}{0.5,0,0.5}

\def\pier #1{{\color{red}#1}}
\def\juerg #1{{\color{red}#1}}
\def\gianni #1{{\color{ciclamino}#1}}
\def\juergen #1{{\color{blue}#1}}

\def\pcol #1{{\color{ciclamino}#1}}

\def\pcol #1{#1}
\def\juergen #1{#1}
\def\pier #1{#1}
\def\juerg #1{#1}
\def\gianni #1{#1}




\bibliographystyle{plain}


%

\finqui

\def\Beq{\Begin{equation}}
\def\Eeq{\End{equation}}
\def\Bsist{\Begin{eqnarray}}
\def\Esist{\End{eqnarray}}

\def\Bthm{\Begin{theorem}}
\def\Ethm{\End{theorem}}
\def\Blem{\Begin{lemma}}
\def\Elem{\End{lemma}}

\def\Brem{\Begin{remark}\rm}
\def\Erem{\End{remark}}

\def\Bdim{\Begin{proof}}
\def\Edim{\End{proof}}
\def\Bcenter{\Begin{center}}
\def\Ecenter{\End{center}}
\let\non\nonumber




\def\step #1 \par{\medskip\noindent{\bf #1.}\quad}


\def\Lip{Lip\-schitz}

\def\aand{\quad\hbox{and}\quad}

\def\lhs{left-hand side}
\def\rhs{right-hand side}



\def\multibold #1{\def\arg{#1}%
  \ifx\arg\pto \let\next\relax
  \else
  \def\next{\expandafter
    \def\csname #1#1#1\endcsname{{\bf #1}}%
    \multibold}%
  \fi \next}

\def\pto{.}

\def\multical #1{\def\arg{#1}%
  \ifx\arg\pto \let\next\relax
  \else
  \def\next{\expandafter
    \def\csname cal#1\endcsname{{\cal #1}}%
    \multical}%
  \fi \next}


\def\multimathop #1 {\def\arg{#1}%
  \ifx\arg\pto \let\next\relax
  \else
  \def\next{\expandafter
    \def\csname #1\endcsname{\mathop{\rm #1}\nolimits}%
    \multimathop}%
  \fi \next}

\multibold
qwertyuiopasdfghjklzxcvbnmQWERTYUIOPASDFGHJKLZXCVBNM.

\multical
QWERTYUIOPASDFGHJKLZXCVBNM.

\multimathop
diag dist div dom mean meas sign supp .


\def\accorpa #1#2{\eqref{#1}--\eqref{#2}}
\def\Accorpa #1#2 #3 {\gdef #1{\eqref{#2}--\eqref{#3}}%
  \wlog{}\wlog{\string #1 -> #2 - #3}\wlog{}}


\def\separa{\noalign{\allowbreak}}

\def\somma #1#2#3{\sum_{#1=#2}^{#3}}
\def\tonde #1{\left(#1\right)}

\def\<#1>{\mathopen\langle #1\mathclose\rangle}
\def\norma #1{\mathopen \| #1\mathclose \|}

\def\[#1]{\mathopen\langle\!\langle #1\mathclose\rangle\!\rangle}

\def\iot {\int_0^t}
\def\ioT {\int_0^T}
\def\intQt{\int_{Q_t}}
\def\intQ{\int_Q}
\def\iO{\int_\Omega}

\def\dt{\partial_t}
\def\dn{\partial_{\mathbf{n}}}
\def\ds{\,ds}

\def\cpto{\,\cdot\,}

\def\checkmmode #1{\relax\ifmmode\hbox{#1}\else{#1}\fi}
\def\aeO{\checkmmode{a.e.\ in~$\Omega$}}
\def\aeQ{\checkmmode{a.e.\ in~$Q$}}

\def\aat{\checkmmode{for a.e.~$t\in(0,T)$}}


\def\erre{{\mathbb{R}}}

\def\enne{{\mathbb{N}}}




\def\genspazio #1#2#3#4#5{#1^{#2}(#5,#4;#3)}
\def\spazio #1#2#3{\genspazio {#1}{#2}{#3}T0}

\def\L {\spazio L}
\def\H {\spazio H}


\def\Lx #1{L^{#1}(\Omega)}
\def\Hx #1{H^{#1}(\Omega)}

\def\LQ #1{L^{#1}(Q)}

\def\Luno{\Lx 1}
\def\Ldue{\Lx 2}

\def\Huno{\Hx 1}
\def\Hdue{\Hx 2}
\def\Hunoz{{H^1_0(\Omega)}}


\def\LQ #1{L^{#1}(Q)}


\let\badtheta\theta
\let\theta\vartheta
\let\eps\varepsilon
\let\phi\varphi

\let\hat\widehat

\let\TeXchi\chi                         
\newbox\chibox
\setbox0 \hbox{\mathsurround0pt $\TeXchi$}
\setbox\chibox \hbox{\raise\dp0 \box 0 }
\def\chi{\copy\chibox}


\def\pn{\mathbf{n}}

\def\VA #1{V_A^{#1}}
\def\VB #1{V_B^{#1}}
\def\VC #1{V_C^{#1}}

\def\Ar{{A^\rho}}
\def\Azr{{A^{2\rho}}}
\def\Bs{{B^\sigma}}
\def\Bzs{{B^{2\sigma}}}
\def\Ct{{C^\tau}}

\def\VAr{{V_A^{\rho}}}

\def\VBs{{V_B^\sigma}}

\def\VCt{{V_C^\tau}}

\def\ful{f_1^\lambda}
\def\Ful{F_1^\lambda}
\def\Fl{F^\lambda}
\def\fl{f^\lambda}

\def\phil{\phi^\lambda}
\def\mul{\mu^\lambda}
\def\Sl{S^\lambda}

\def\phiab{\phi_{\alpha,\beta}}
\def\muab{\mu_{\alpha,\beta}}
\def\Sab{S_{\alpha,\beta}}

\def\phiz{\phi_0}
\def\muz{\mu_0}
\def\Sz{S_0}

\def\soluz{(\mu,\phi,S)}

\Begin{document}


%
\title{Asymptotic analysis of a tumor growth model\\ with fractional operators}
\author{}
\date{}

\maketitle
\Bcenter
\vskip-1cm
{\large\sc Pierluigi Colli$^{(1)}$}\\
{\normalsize e-mail: {\tt pierluigi.colli@unipv.it}}\\[.25cm]
{\large\sc Gianni Gilardi$^{(1)}$}\\
{\normalsize e-mail: {\tt gianni.gilardi@unipv.it}}\\[.25cm]
{\large\sc J\"urgen Sprekels$^{(2)}$}\\
{\normalsize e-mail: {\tt sprekels@wias-berlin.de}}\\[.45cm]
$^{(1)}$
{\small Dipartimento di Matematica ``F. Casorati'', Universit\`a di Pavia}\\
{\small and Research Associate at the IMATI -- C.N.R. Pavia}\\
{\small via Ferrata 5, 27100 Pavia, Italy}\\[.2cm]
$^{(2)}$
{\small Department of Mathematics}\\
{\small Humboldt-Universit\"at zu Berlin}\\
{\small Unter den Linden 6, 10099 Berlin, Germany}\\[2mm]
{\small and}\\[2mm]
{\small Weierstrass Institute for Applied Analysis and Stochastics}\\
{\small Mohrenstrasse 39, 10117 Berlin, Germany}
\Ecenter

%
\Begin{abstract}\noindent
In this paper, we study a system of three evolutionary operator equations involving fractional powers
of selfadjoint, monotone, unbounded, linear operators having compact resolvents. 
This system constitutes a generalized and relaxed  
version of a phase field system of Cahn--Hilliard type modelling tumor growth that has 
originally been proposed 
in Hawkins-\juergen{Daarud} et al. ({\em Int.\,\,J. Numer.\,\,Math.\,\,Biomed.\,\,Eng.\,\,{\bf 28} (2012), 3--24}). The original phase field system and certain relaxed versions thereof have been 
studied in recent papers co-authored by the present authors and E. Rocca. The model consists of a 
Cahn--Hilliard equation for the tumor cell fraction \,$\varphi$, coupled to a reaction-diffusion equation for a function $\,S\,$ representing the nutrient-rich extracellular water volume fraction. Effects due to fluid motion are neglected. Motivated by the possibility that the diffusional regimes governing the evolution of the different constituents of the model may be of different (e.g., fractional) type, the present authors studied in a recent note a generalization of the 
systems investigated in the abovementioned works. Under rather general assumptions, well-posedness and regularity results have been shown. In particular, by 
writing the equation governing the evolution of the chemical potential in the form of a general variational 
inequality, also singular or nonsmooth contributions of
logarithmic or of double obstacle type to the energy density could be admitted. In this note, we perform
an asymptotic analysis of the governing system as two (small) relaxation parameters approach zero separately and simultaneously. Corresponding well-posedness and regularity results are 
established for the respective cases; in particular, we give a detailed discussion which assumptions on the
admissible nonlinearities have to be postulated in each of the occurring cases.  
\vskip3mm
\noindent {\bf Key words:}
Fractional operators, Cahn--Hilliard systems, well-posedness, regularity
of solutions, tumor growth models, asymptotic analysis. 
\vskip3mm
\noindent {\bf AMS (MOS) Subject Classification:} {35B40, \pier{35K55}, 35K90, 35Q92, 92C17.}
\End{abstract}
\salta
\pagestyle{myheadings}
\newcommand\testopari{\sc Colli \ --- \ Gilardi \ --- \ Sprekels}
\newcommand\testodispari{\sc Asymptotic analysis of a  fractional tumor growth model}
\markboth{\testopari}{\testodispari}
\finqui
%

\section{Introduction}
\label{Intro}
\setcounter{equation}{0}
Let $\Omega\subset\erre^3$ denote an open, bounded, and connected set
with smooth boundary $\Gamma$ and unit outward normal $\pn$\pier{;} let $T>0$ be given. 
\pier{Setting} 
$Q_t:=\Omega\times (0,t)$ for $t\in (0,T)$ and  $Q:=\Omega\times (0,T)$, as well as 
$\Sigma:=\Gamma\times (0,T)$\pier{, we} investigate in this paper the evolutionary system
\Bsist
  && \alpha\, \dt\mu + \dt\phi + A^{2\rho} \mu 
  = P(\phi)(S-\mu)\qquad\mbox{in \,$Q$,}
  \label{Iprima}
  \\[1mm]
  && \mu =  \beta\, \dt\phi + B^{2\sigma}\phi + F'(\phi) \qquad\mbox{in \,$Q$},
  \label{Iseconda}
  \\[1mm]
  && \dt S + C^{2\tau}S 
  = - P(\phi)(S-\mu)\qquad\mbox{in \,$Q$,}
  \label{Iterza}
  \\[1mm]
  && \mu(0) = \muz, \quad
  \phi(0) = \phiz, \quad
  S(0) = \Sz,\qquad\mbox{in \,$\Omega$.}
  \label{Icauchy}
\Esist

In the above system,  $A^{2\rho}$, $B^{2\sigma}$, $C^{2\tau}$, with $\rho>0,\,\sigma>0,\,\tau>0$,
denote fractional powers of the selfadjoint,
monotone, and unbounded linear operators $A$, $B$, and~$C$, respectively, 
which are supposed to be densely defined in $H:=\Ldue$ and to have compact resolvents.
Moreover, $F'$~denotes the derivative of a double-well potential~$F$.
Typical and physically significant examples of $F$ 
are the so-called {\em classical regular potential}, the {\em logarithmic double-well potential\/},
and the {\em double obstacle potential\/}, which are given, in this order,~by
\begin{align}
\label{regpot}
  & F_{reg}(r) := \frac 14 \, (r^2-1)^2 \,,
  \quad r \in \erre,\\[1mm] 
  \label{logpot}
  & F_{log}(r) := \left\{\begin{array}{ll}
  (1+r)\ln (1+r)+(1-r)\ln (1-r) - c_1 r^2\,,&
  \quad r \in (-1,1)\\
  2\log(2)-c_1\,,&\quad r\in\{-1,1\}\\
  +\infty\,,&\quad r\not\in [-1,1]
\end{array}\right. \,,
  \\[1mm]
\label{obspot}
  & F_{2obs}(r) := \gianni{c_2 (1- r^2)}
  \quad \hbox{if $|r|\leq1$}
  \aand
  F_{2obs}(r) := +\infty
  \quad \hbox{if $|r|>1$}.
\end{align}
Here, the constants $c_i$ in \eqref{logpot} and \eqref{obspot} satisfy
$c_1>1$ and $c_2>0$, so that the corresponding functions are nonconvex.
In cases like \eqref{obspot}, one has to split $F$ into a nondifferentiable convex part~$F_1$ 
(the~indicator function of $[-1,1]$, in the present example) and a smooth perturbation~$F_2$.
Accordingly, in the term $F'(\phi)$ appearing in~\eqref{Iseconda}, 
one has to replace the derivative $F_1'$ of the convex part $F_1$
by the subdifferential $\partial F_1$ and interpret \eqref{Iseconda} as a differential inclusion
or as a \juergen{variational} inequality involving $F_1$ rather than~$F_1'$.
Furthermore, the function $P$ occurring in \eqref{Iprima} and \eqref{Iterza}
is nonnegative and smooth, and the terms on the \rhs s in \eqref{Icauchy} are prescribed initial data.

The above system is a generalization of a system of PDEs that \pier{constitutes} a relaxed version of a model
for tumor growth originally introduced in \cite{HZO12} that was investigated in the 
papers\pier{\cite{CGRS, CGRS1, CGRS2}} co-authored by the authors of this note and E. Rocca. In these works,
we studied the special situation when $A^{2\rho}=B^{2\sigma}=C^{2\tau}=-\Delta$ with zero 
Neumann boundary conditions, and established general results concerning well-posedness,
regularity, and optimal control. In particular, in\pier{\cite{CGRS,CGRS2}} a thorough asymptotic analysis,
coupled with rigorous error estimates,
was performed for the situation when the relaxation parameters $\alpha>0$ and $\beta>0$
approach zero, either separately or simultaneously. Notice also that in the case $P\equiv 0$ the equation
\eqref{Iterza} decouples from the other two equations \eqref{Iprima}, \eqref{Iseconda};
the latter system of equations has for the case $\alpha=0$ recently been the subject of a series of
investigations by the present authors (cf. the papers \cite{CGS18,CGS19,CGS21,CGS22}).   
 
In this paper, we intend to perform a
corresponding asymptotic analysis for the general system \eqref{Iprima}--\eqref{Icauchy}, where we take
advantage of the well-posedness and regularity results that were established in \pier{our recent paper \cite{CGS23}}. It  will be demonstrated that for each of the three limit processes 
$$\alpha\searrow0\,,\,\beta>0,\qquad \alpha>0\,,\,\beta\searrow0,\qquad \alpha\searrow0\,,\beta\searrow0,$$
meaningful limit problems \pier{occur} for which the existence of solutions can be shown. In this analysis, it will
turn out that each of the three limit processes needs specific assumptions for the fractional operators and the
admissible nonlinearities. We will also address questions of uniqueness and continuous dependence, where, again,
specific assumptions are necessary for the three cases.

Modeling the dynamics of tumor growth has recently become an important issue in applied mathematics (see, e.g.,~\cite{CL2010,Lowen08}), and \pier{some different models have been 
introduced and discussed, numerical simulations have been provided and a comparison 
with the behavior of other special materials has been in order; for all that we just refer to, e.g., the works \cite{BLM,CLLW,CL2010,FBG2006,Fri2007,Lowen10,HDPZO,OHP,WZZ}.
In particular, about diffuse interface models, we point out that these models mostly 
follow the Cahn--Hilliard framework (see \cite{Cahn}) that originated 
from the theory of phase 
transitions and is extensively employed in 
materials science and multiphase fluid flow.
Among these models, two main classes can be categorized: 
the first one looks at the tumor and healthy cells as inertialess fluids and 
takes the effects generated by the fluid flow development into account by  
postulating a Darcy or a Brinkman law; in this direction, we refer to
\juergen{\cite{DFRSS,EGAR,FLRS,GL2016,GL2018,GLNS, GLSS, JWZ,LTZ,SW, Wise2011}} 
(cf.~\cite{BCG,ConGio,DG,DGG,FW2012,GioGrWu,WW2012, WZ2013} as well, where local or nonlocal Cahn--Hilliard systems with Darcy or Brinkman law are dealt with). Moreover,
further mechanisms such as chemotaxis and active transport can be considered in the phenomenology. On the other hand, the second class of models, including the one 
from which the system \eqref{Iprima}--\eqref{Icauchy} originates, actually neglects the velocity and admits as variables concentrations and chemical potential.  
A variety of contributions inside this class is provided by the works\cite{CRW,CWSL,CGMR, FLR, GL2017-1, GL2017-2,GLR, MRS, S, S_DQ,S_b,S_a}.}

\pier{To our knowledge, except for the recent papers\pier{\cite{CGS23,CGS25}}, 
fractional operators have not been 
studied in either of these two groups of models, although 
one may also wonder about nonlocal operators. We point out that 
in recent years fractional operators provide a challenging subject for 
mathematicians: they have been successfully utilized in many different situations, and 
a wide literature already exists about equations and systems with fractional terms.
For an overview of recent contributions, we refer to the papers 
\cite{CGS18, CGS19} and \cite{CG1}, which offer to the interested reader
a number of suggestions to deepen the knowledge of the field. 
In our approach here, we adopt the setting 
of~\cite{CGS23} and consequently work on fractional 
operators defined via spectral theory. This framework includes, in particular, 
powers of a second-order elliptic operator with either Dirichlet or Neumann or 
Robin \juergen{homogeneous} boundary conditions, and other operators like, e.g., 
fourth-order ones or systems involving the Stokes operator. 
A precise definition for our fractional operators
$A^{2\rho}$, $B^{2\sigma}$, $C^{2\tau}$ along with their properties  
will be  given in the first part of Section~2 below.} 
\pcol{As far as a biological background 
for the system \accorpa{Iprima}{Iterza} is concerned, we 
claim that in our approach the three fractional operators, 
which may be considerably 
different the one from the other, are employed for the dynamics 
of tumor growth and diffusion processes. 
The three operators $A^{2\rho}$, $B^{2\sigma}$, $C^{2\tau}$ 
are allowed to be a variation of fractional Laplacians, 
but also other elliptic operators, and may show different 
orders.  Indeed, some components 
in tumor development, such as immune cells, exhibit an anomalous 
diffusion dynamics (as it observed in experiments~\cite{EGPS}), but other 
components, like chemical potential and nutrient concentration 
are possibly governed by different fractional or non-fractional flows. 
However, taking all this into account, it is the case of pointing out that 
fractional operators are becoming more and more implemented in the field of
biological applications: to this concern, a selection of notable 
and meaningful references is given by \cite{BKM, CB, EGPS, EL, G-B, INK, JJ, KU, 
MV, SAJM, SA, ZSMD}.}

The paper is organized as follows: 
in the next section, we list our assumptions and notations, and we state some results
for the system \eqref{Iprima}--\eqref{Icauchy} that \pier{are} valid if both 
\juergen{of} the relaxation parameters are positive. The following
sections then bring the asymptotic analysis as the parameters $\alpha>0$ and $\beta>0$ approach
zero, where each of the relevant cases will be treated in a \juergen{separate} section.
Throughout this paper, we make use of
the elementary Young inequality
\begin{equation}
\label{young}
 ab\leq \gamma a^2 + \frac 1 {4\gamma}\,b^2
  \quad \hbox{for every $a,b\in\erre$ and $\gamma>0$}.
\end{equation}
Moreover, given a Banach space $\,X$, we denote by $\|\,\cdot\,\|_X\,$ its norm and by $\,X^*\,$ its dual.
The dual pairing between $X^*$ and $X$ is denoted by $\,\langle\,\cdot\,,\,\cdot\,\rangle_X$. 
The only exception from this rule is the space $H := \Ldue$, for which 
 $\|\,\cdot\,\|$ and $(\,\cdot\,,\,\cdot\,)$ denote the standard norm and inner product, respectively.

\section{General assumptions and known results}
\setcounter{equation}{0}

In this section, we give precise assumptions and notations and state some  results for
the relaxed system where $\alpha>0$ and $\beta>0$. 
Now, we start introducing our assumptions.
As for the operators, we first postulate that

\vspace{3mm}\noindent
{\bf (A1)} \,\,\,$A:D(A)\subset H\to H$, $B:D(B)\subset H\to H$, and $C:D(C)\subset H\to H$,
 are \hspace*{13.5mm}unbounded, monotone, selfadjoint, linear operators with compact resolvents.
  
\vspace{3mm}\noindent	
Therefore, there are sequences 
$\{\lambda_j\}$, $\{\lambda'_j\}$, $\{\lambda''_j\}$, and $\{e_j\}$, $\{e'_j\}$, $\{e''_j\}$, 
of eigenvalues and corresponding eigenfunctions such that
\begin{align}
\label{eigen}
& A e_j = \lambda_j e_j, \quad
 B e'_j = \lambda'_j e'_j,
\quad C e''_j = \lambda''_j e''_j, 
\quad \hbox{with} \quad
(e_i,e_j) = (e'_i,e'_j) = (e''_i,e''_j) = \delta_{ij},\nonumber\\
&\qquad\hbox{for  all $\,i,j\in \enne$}\,,  
\\[1mm]
\separa
\label{eigenvalues}
& 0 \leq \lambda_1 \leq \lambda_2 \leq \dots , \quad
0 \leq \lambda'_1 \leq \lambda'_2 \leq \dots,
\aand
0 \leq \lambda''_1 \leq \lambda''_2 \leq \dots , \quad\hbox{where} 
\non
\\
& \qquad 
\lim_{j\to\infty} \lambda_j
= \lim_{j\to\infty} \lambda'_j
= \lim_{j\to\infty} \lambda''_j
= + \infty\,,
  \\[1mm]
  \separa
\label{complete}
  & \hbox{$\{e_j\}$, $\{e'_j\}$, and $\{e''_j\}$, are complete systems in $H$}.
\end{align}
As a consequence, we can define the powers of these operators with arbitrary 
positive real exponents as done below.
As far as the first operator is concerned, we have for $\rho>0$
\Bsist
  && \VAr := D(\Ar)
  = \Bigl\{ v\in H:\ \somma j1\infty |\lambda_j^\rho (v,e_j)|^2 < +\infty \Bigr\}
  \aand
  \label{defdomAr}
  \\[-3mm]
  && \Ar v = \somma j1\infty \lambda_j^\rho (v,e_j) e_j
  \quad \hbox{for $v\in\VAr$},
  \label{defAr}
\Esist
the series being convergent in the strong topology of~$H$,
due to the properties \eqref{defdomAr} of the coefficients.
We endow $\VAr$ with the graph norm, i.e., we~set
\Beq
  (v,w)_\VAr := (v,w) + (\Ar v,\Ar w)
  \aand
  \norma v_\VAr := (v,v)_\VAr^{1/2}
  \quad \hbox{for $v,w\in\VAr$},
  \label{defnormaAr}
\Eeq
and obtain a Hilbert space.
In the same way, 
we can define the powers $B^\sigma$ and $C^\tau$ for every $\sigma>0$ and $\tau>0$,
starting from \eqref{eigen}--\eqref{complete} for $B$ and~$C$.
We therefore set
\,$\VBs := D(\Bs)$ and $\VCt := D(\Ct)$, endowed with the norms 
$\|\cdot\|_\VBs$ and $\|\cdot\|_\VCt$ induced by the inner products
\begin{align}
   \label{defprodBsCr}\non
   & (v,w)_\VBs := (v,w) + (B^\sigma v,B^\sigma w)
  \aand
  (v,w)_\VCt := (v,w) + (C^\tau v,C^\tau w),
  \non
  \\
  & \quad \hbox{for $v,w\in \VBs$ and $v,w\in\VCt$, respectively}.
 \end{align}
Since $\lambda_j\geq0$ for every~$j$, one immediately deduces from the definition of \pier{$A^\rho$} that
\Bsist
  && \Ar : \VAr \subset H \to H
  \quad \hbox{is maximal monotone, and} 
  \non
  \\
  &&{\eps I + \Ar} : \VAr \to H
  \quad \hbox{is for every $\,\eps>0\,$ a topological isomorphism with the inverse}\nonumber\\
  &&(\eps I+\Ar)^{-1}v=\sum_{j=0}^\infty \bigl(\eps+ \lambda_j^r\bigr)^{-1}(v,e_j)e_j\quad\mbox{for }\,v\in H,
  \label{maxmon}
\Esist
where $I:H\to H$ is the identity operator. Similar results hold 
for $B^\sigma$ and~$C^\tau$.
It is clear that, for every $\rho_1,\,\rho_2>0$, we have the Green type formula
\Beq
  (A^{\rho_1+\rho_2} v,w)
  = (A^{\rho_1} v, A^{\rho_2} w)
  \quad \hbox{for every $v\in\VA{\rho_1+\rho_2}$ and $w\in V_A^{\rho_2}$},
  \label{propA}
\Eeq
and that similar relations holds for the other two types of fractional operators.
Due to these properties, we can define proper extensions of the operators that
allow values in dual spaces.
In particular, we can write variational formulations of the equation \accorpa{Iprima}{Iterza}.
It is convenient to use the notations
\Beq
  V_A^{-\rho} := (\VAr)^* , \quad
  V_B^{-\sigma} := (\VBs)^*, 
  \aand
  V_C^{-\tau} := (\VCt)^* ,
  \quad \hbox{for $\rho>0,\,\sigma>0,\,\tau>0$}.
  \label{defnegspaces}
\Eeq
Thus, we have that
\Beq
  \Azr \in \calL(\VAr, V_A^{-\rho}) , \quad
  \Bzs \in \calL(\VBs, V_B^{-\sigma}),
  \aand
  C^{2\tau} \in \calL(\VCt, V_C^{-\tau}) ,
  \label{extensions2}
\Eeq
as well as
\Beq
 \Ar \in \calL(H,V_A^{-\rho}) , \quad
  B^\sigma \in \calL(H, V_B^{-\sigma}),
  \aand
  C^\tau \in \calL(H, V_C^{-\tau}) .
  \label{extensions1}
\Eeq
The symbols $\<\cpto,\cpto>_\VAr$ and $\,\langle\,\cdot\,,\,\cdot\,\rangle_{\VCt}$ will be used for the duality pairings
between $V_A^{-\rho}$ and~$\VAr$ and between $V_C^{-\tau}$ and $\VCt$, respectively. 
Moreover, we identify $H$ with a subspace of $V_A^{-\rho}$
in the usual way, i.e., such that
\Beq
  \< v,w >_\VAr = (v,w)
  \quad \hbox{for every $v\in H$ and $w\in\VAr$}.
  \label{identification}
\Eeq
Analogously, we have that $H\subset V_B^{-\sigma}$ and $H\subset V_C^{-\rho}$ and use similar notations.
Notice (see, e.g., \cite[Sect.~3]{CGS18}) that all of the embeddings
\begin{align}
\label{embedA}
  & \VA{r_2} \subset \VA{r_1} \subset H, \quad\mbox{for }\,0<r_1<r_2,\\
  \label{embedAm}
    & H \subset \VA{-r_1} \subset \VA{-r_2},\quad\mbox{for }\,0<r_1<r_2,\\
  \label{embedB}
  & \VB{\sigma_2} \subset \VB{\sigma_1} \subset H,\quad\mbox{for } \,0<\sigma_1<\sigma_2, \\
\label{embedC}
&\juergen{V^{\tau_2}_C\subset V^{\tau_1}_C \subset H,\quad\mbox{for }\,0<\tau_1<\tau_2,}
\end{align}
are dense and compact. 

\vspace{3mm}
\gianni{From now on, we assume:}

\vspace{3mm} \noindent
{\bf (A2)} \,\,\, $\rho$, $\sigma$ and $\tau$ are fixed positive real numbers.

\vspace{3mm}\noindent
For the  nonlinear functions entering the equations \accorpa{Iprima}{Iterza} of our system, we postulate the properties listed below:

\vspace{3mm}\noindent
{\bf (A3)} \,\,\,$ F = F_1 + F_2,$ where:
  \begin{align}
  \label{hpFuno}
  & F_1 : \erre \to [0,+\infty]
  \quad \hbox{is convex and lower semicontinuous with} \quad
  F_1(0) = 0.
    \\[1mm]
  \label{hpFdue}
  & F_2 \in C^1(\erre), \,\,  \hbox{ and }\,\,F_2' \,\,\hbox{ is Lipschitz continuous with 
  Lipschitz constant }\, L>0.
  		\\[1mm]
\label{hpbelow}
  & \gianni{F(s) \geq \overline c_1 s^2 - \overline c_2 \quad
  \hbox{for some positive constants $\overline c_1$ and $\overline c_2$ and every $s\in\erre$}.}
    	\\[1mm]
 \label{hpP} 
	& P :\erre \to [0,+\infty) \quad \hbox{is bounded and \Lip\ continuous}.
\end{align}

\vspace{2mm}
We set, for convenience,
\Beq
\label{deffunodue}  
  f_1 := \partial F_1 
  \aand
 f_2 := F_2' \,,
\Eeq
and denote by $D(F_1)$ and $D(f_1)$ the effective domain of $F_1$ and~$f_1$, respectively.
We notice that $f_1$ is a maximal monotone graph in $\erre\times\erre$
and use the same symbol $f_1$ for the maximal monotone operators induced in $L^2$ spaces. For
every \,$s\in D(f_1)$, we denote by $\,f_1^\circ(s)\,$ the element of minimal modulus
in $\,f_1(s)$. Moreover, if
 the subdifferential $\,\partial F_1(s)\,$ is a singleton for every $s\in D(f_1)$ 
 (which is, e.g., the case if $F_1\in C^1(\erre)$),
then we identify the singleton $\{f_1(s)\}$ with the real number $f_1(s)$ and treat 
the mapping $s\mapsto f_1(s)$ as a real-valued function without further comment. 

\vspace{3mm} 
Using \eqref{propA} and its analogues for~$B$ and~$C$, we can give a weak formulation 
of the equations \accorpa{Iprima}{Iterza}.
Moreover, we present \eqref{Iseconda} as a variational inequality.
For the data, we make the following assumptions:

\vspace{3mm}\noindent
{\bf (A4)} \,\,\,$\muz \in H$, \,\,\,
  $\phiz \in \VBs$   \,\,\,with\,\,\,
  $F_1(\phiz) \in \Luno $, 
  \,\, and\, $\Sz \in H$ .

\vspace{3mm}\noindent
\gianni{By assuming $\alpha\geq0$ and $\beta\geq0$, we then look for a triple $\soluz$ satisfying
\begin{align}
\label{regmu} 
&\mu\in \gianni{\L2{\VAr}},\\[1mm]
\label{regphi}  
& \phi \in  \L\infty{\VBs}, \quad
 \beta\,\dt\phi\in L^2(0,T;H),
\\[1mm]
\label{regdtmp}
& \dt(\alpha\mu+\phi) \in L^2(0,T;V_A^{-\rho}),\\[1mm]
  \label{regS}
& S \in \H1{V_C^{-\tau}} \cap \L\infty H \cap \L2{\VCt},\\[1mm]
  \label{L1Q}
& F_1(\phi) \in \LQ1,
\end{align} 
}%
and solving the system
\begin{align}
\label{prima}
&\< \,\dt(\alpha\,\mu(t)+\phi(t)) , v\, >_\VAr
  \,+ \, ( \Ar \mu(t) , \Ar v )
  \,=\, \bigl( P(\phi(t)) ( S(t) - \mu(t) ) , v \bigr)
  \nonumber\\[1mm]
& \qquad \hbox{for every $v\in\VAr$ and \aat,}\\[2mm]
  \label{seconda}
& \bigl( \beta\,\dt\phi(t) , \phi(t) - v \bigr)
  + \bigl( B^\sigma \phi(t) , B^\sigma( \phi(t)-v) \bigr)
  \nonumber   \\[1mm]
  & 
  + \iO F_1(\phi(t))
  + \bigl( f_2(\phi(t)) ,  \phi(t)-v \bigr)
  \,\leq\, \bigl( \mu(t) , \phi(t)-v \bigr)
  + \iO F_1(v)
  \nonumber\\[1mm]
  & \qquad \hbox{for every $v\in\VBs$ and \aat,}\\[2mm]
  \label{terza}
& \< \,\dt S(t) , v \,>_\VCt
  \,+\, ( \Ct S(t) ,C^\tau v )
  \,=\, - \bigl( P(\phi(t)) ( S(t) - \mu(t) ) , v \bigr)
  \nonumber
  \\[1mm]
  & \qquad \hbox{for every $v\in\VCt$ and \aat},\\[2mm]
  \label{cauchy}
& \gianni{(\alpha\mu+\phi)(0) = \alpha\muz+\phiz \,, \quad
  (\beta\phi)(0) = \beta\phiz\,,} 
  \aand
  S(0) = \Sz \,.
\end{align}
Here, it is understood that $\,\,\iO F_1(v)=+\infty\,$ whenever $\,F_1(v)\not\in
L^1(\Omega)$.

\gianni{%
\Brem
\label{Remcauchy}
The above formulation is meaningful for nonnegative coefficients $\alpha$ and~$\beta$.
This holds, in particular, for~\eqref{cauchy}. 
However, \juergen{depending} on whether these coefficients are positive or zero, 
the initial conditions can be reformulated in a more explicit way, namely,
\begin{align}
  & \mu(0) = \muz , \quad
  \phi(0) = \phiz,
  \aand
  S(0) = \Sz,
  && \hbox{if $\alpha>0$ and $\beta>0$},
  \label{cauchyabpos}
  \\
  & (\alpha\mu+\phi)(0) = \alpha\muz+\phiz,
  \aand
  S(0) = \Sz,
  && \hbox{if $\alpha>0$ and $\beta=0$},
  \label{cauchyaposbzero}
  \\
  & \phi(0) = \phiz,
  \aand
  S(0) = \Sz,
  && \hbox{if $\alpha=0$ and $\beta\geq 0$}.
  \label{cauchyazero}
\end{align}
\Erem
}%

\vspace{2mm}\noindent
Observe that \eqref{prima}--\eqref{terza} are equivalent to their time-integrated variants, \pier{in particular for \eqref{seconda} we have}
\Bsist
  && \ioT \bigl( \beta\,\dt\phi(t) , \phi(t) - v(t) \bigr) \, dt 
  + \ioT \bigl( B^\sigma\phi(t) , B^\sigma(\phi(t)-v(t)) \bigr)\,dt 
  \nonumber
  \\ 
  && 
  + \intQ F_1(\phi) 
  + \ioT \bigl( f_2(\phi(t)) , \phi(t)-v(t) \bigr) \, dt 
  \non
  \\
  && \leq \ioT \bigl( \mu(t) , \phi(t)-v(t) \bigr)\,dt
  + \intQ F_1(v)
  \quad \hbox{for every $v\in\L2{\VB\sigma}$},
  \qquad
  \label{intseconda}
\Esist
where we put \,$\int_Q F_1(v)=+\infty\,$ whenever $\,F_1(v)\not\in L^1(Q)$.

The following result was proved in \cite[Thms.~2.3~and~2.5]{CGS23}:

\Bthm
\label{Wellposed}
Let the assumptions {\bf (A1)}--{\bf (A4)} be fulfilled,
and assume that $\alpha>0$ and $\beta>0$.
Then there exists a triple $(\mu,\phi,S)$ with the regularity \juergen{\eqref{regmu}}--\eqref{L1Q} that solves the problem 
\gianni{\eqref{prima}--\eqref{terza} and the initial conditions \eqref{cauchyabpos}}.
Moreover, this solution satisfies the estimate
\begin{align} 
\label{bounds1}
& \|\dt(\alpha\,\mu+\phi)\|_{L^2(0,T;V_A^{-\rho})} \,+\, \alpha^{1/2} \norma\mu_{\L\infty H}
  \,+\, \norma{\Ar\mu}_{\L2H}
  \nonumber\\[1mm]
&   + \,\beta^{1/2}\, \norma{\dt\phi}_{\L2H}
  \,+\, \|\phi\|_{L^\infty(0,T;\VBs)}
 \, +\, \gianni{\norma{F(\phi)}_{\L\infty\Luno}}
  \nonumber\\[1mm]
&   +\, \norma S_{H^1(0,T;V_C^{-\tau})\cap C^0([0,T];H)\cap L^2(0,T;\VCt)}
    \,+\, \norma{P^{1/2}(\phi)(S-\mu)}_{\L2H}
  \nonumber
  \\[1mm]
  & \leq \,\widehat K_1 \,\bigl(
    \alpha^{1/2} \norma\muz
    \,+\, \norma{B^\sigma\phiz}
    \,+\, \norma{F(\phiz)}_{\Luno}
    \,+\, \norma\Sz
    + 1  \bigr)\,,
\end{align}
with a constant~$\,\widehat K_1>0\,$  that depends only on~$\Omega$,
\gianni{the constants $\overline c_1$ and $\overline c_2$ from \eqref{hpbelow}, and~$P$}. 
If, in addition, the condition
\begin{equation}
\label{ini2}
\mu_0\in\VAr,\quad \phiz\in V_B^{2\sigma}\,\,\mbox{ with }\,\,f_1^\circ(\phiz)\in H,\quad
S_0\in\VCt,
\end{equation}
is fulfilled, \juergen{then} the above solution enjoys the further regularity
\begin{align}
\label{adregmu}
&\mu\in H^1(0,T;H)\cap L^\infty(0,T;\VAr)\cap L^2(0,T;V_A^{2\rho}),\\
\label{adregphi}
&\phi\in W^{1,\infty}(0,T;H)\cap H^1(0,T;\VBs),\\
\label{adregS}
&S\in H^1(0,T;H)\cap L^\infty(0,T;\VCt)\cap L^2(0,T;V_C^{2\tau}).
\end{align}
Moreover, if the embedding conditions 
\Beq
  \gianni{\VA\rho} \subset \Lx 4
  \aand
  \gianni{\VC\tau} \in \Lx 4 \,.
  \label{embedAC}
\Eeq
are fulfilled, then the above solution is uniquely determined. 
\Ethm

\vspace{2mm}
\Brem
\label{RemembedAC}
The first embedding in \eqref{embedAC} is, for instance,  satisfied 
if $\Azr \pier{{}= A:{}} =-\Delta$ with the domain $\Hdue\cap\Hunoz$ 
(thus, with zero Dirichlet conditions, but similarly for 
zero Neumann boundary conditions \pier{with domain $\{v\in\Hdue : \ \dn v =0\, $ on 
$\, \Gamma \}$}).
Indeed, we have $\VAr=\Hunoz$ in this case.
Clearly, the same embedding holds true if $\rho$ is sufficiently close to~$\pier{1/2}$.
\Erem

\Brem
\label{Controls}
More generally, we could add known forcing terms $u_\mu$, $u_\phi$ and $u_S$
to the \rhs s of equations \eqref{Iprima}, \eqref{Iseconda} and~\eqref{Iterza}, respectively,
and accordingly modify the definition of solution.
If we assume that
\Beq
  u_\mu ,\ u_\phi \,, u_S \in \L2H\,,
  \label{hpcontrols}
\Eeq
then we have a similar well-posedness result.
In estimate~\eqref{bounds1}, one has to modify the \rhs\ by adding 
the norms corresponding to \eqref{hpcontrols} 
(possibly multiplied by negative powers of $\alpha$ and~$\beta$).
This remark is useful for performing a control theory of the above system
with distributed controls.
\Erem

\Brem
\label{lapprox}
We cannot repeat the proof given in \cite{CGS23}, here. We only note for later use that the result is achieved by
approximation using the Moreau--Yosida regularizations 
$\Ful$ and $\ful$ of $F_1$ of $f_1$ at the level $\lambda>0$ introduced in, e.g., \cite[p.~28 and p.~39]{Brezis}. We  set, 
for convenience,
\Beq
  \label{yosida1}
  \Fl := \Ful + F_2
  \aand
  \fl := \ful + f_2 \,.
  \Eeq
Denoting by $J^\lambda:=(\mbox{id}+\lambda\,f_1)^{-1}$ (where $\mbox{id}:\erre\to\erre$ 
is the identity mapping) the resolvent mapping  associated with the maximal monotone 
graph $\,f_1$ for $\lambda>0$, we recall some well-known properties of this regularization, namely,
\begin{align}
  \label{yosida2}
& \Ful(s) = \int_0^s \ful(s') \, ds' \,, \quad
  0 \leq F_1(J^\lambda(s))\le \Ful(s) \leq F_1(s)\,\quad\mbox{for every }\,s\in\erre,\\[1mm]
\label{yosida3}
&   \ful(s) \in f_1(J^\lambda(s))  \, \quad \mbox{for every }\,s\in\erre\,,
\end{align}
and  it follows from \eqref{hpbelow} that there are constants $\hat c_1>0$, 
$\hat c_2>0$, and $\Lambda>0$, such that, for all $\lambda\in (0,
\Lambda)$, we~have
\begin{equation}\label{yosida4}
   \Fl(s) \geq \hat c_1\,s^2\,-\, \widehat c_2
\,\quad\pier{\mbox{for every }}\,s\in\erre\,.
\end{equation}
In the following, we always tacitly assume that $0<\lambda<\Lambda$ when working with
Moreau--Yosida approximations.

Now, we replace $F_1$ in \eqref{seconda} by~$\Ful$ to obtain the system
\Bsist
  && \alpha\, \< \dt\mul(t) , v >_\VAr
  + \bigl( \dt\phil(t) , v \bigr)
  + ( \Ar \mul(t) , \Ar v )
  = \bigl( P(\phil(t)) (\Sl(t)-\mul(t)) , v \bigr)
  \non
  \\
  && \quad \hbox{for every $v\in \VAr$ and \aat},
  \label{primal}
  \\[1mm]
  \separa
  && \beta \,\bigl( \dt\phil(t) , \phil(t) - v \bigr) 
  + \bigl( B^\sigma\phil(t) , B^\sigma(\phil(t)-v) \bigr)
  \non
  \\
  && \quad {}
  + \iO \Ful(\phil(t))
  + \bigl( f_2(\phil(t)) , \phil(t)-v \bigr)
  \non
  \\
  && \leq \bigl( \mul(t) , \phil(t)-v \bigr)
  + \iO \Ful(v)
  \quad \hbox{for every $v\in\VBs$ and \aat},
  \qquad
  \label{secondal}
  \\[1mm]
  \separa
  && \< \dt\Sl(t) , v >_\VCt
  + ( C^\tau\Sl(t) ,C^\tau v )
  = - \bigl( P(\phil(t)) ( \Sl(t) - \mul(t) ) , v \bigr)
  \non
  \\
  && \quad \hbox{for every $v\in\VCt$ and \aat},
  \label{terzal}
  \\[1mm]
  && \mul(0) = \muz \,, \quad
  \phil(0) = \phiz\,, 
  \aand
  \Sl(0) = \Sz \,.
  \label{cauchyl}
\Esist
\Accorpa\Pbll primal cauchyl
Observe that \eqref{secondal} is equivalent to both 
its time-integrated analogue 
and the pointwise variational equation 
(since $\Ful$ is differentiable and $\ful$ is its globally Lipschitz continuous derivative)
\Bsist
  && \beta \bigl( \dt\phil(t) , v \bigr)
  + \bigl( B^\sigma\phil(t) , B^\sigma v \bigr)
  + \bigl( \fl(\phil(t)) ,  v \bigr)
  = \bigl( \mul(t) ,  v \bigr)
  \non
  \\[1mm]
  && \quad \hbox{for every $v\in\VBs$ and \aat}.
  \label{eqsecondal}
\Esist

In the proof of \cite[Thm.~2.3]{CGS23}, it was shown under slightly weaker assumptions on $F$ that 
the system \eqref{primal}, \eqref{terzal}--\eqref{eqsecondal} has 
for every $\lambda\in (0,\Lambda)$ a unique solution triple $(\mul,\phil,\Sl)$ satisfying 
\eqref{regmu}--\eqref{regS} 
and the estimate 
\begin{align}
\label{bounds1l}
&\left\|\dt(\alpha\,\mul+\phil)\right\|_{L^2(0,T;V_A^{-\rho})}\,+\,
\alpha^{1/2}\left\|\mul\right\|_{ L^\infty(0,T;H)}\,+\,\left\|\Ar\mul\right\|_{L^2(0,T;H)}\nonumber\\[1mm]
&+\,\beta^{1/2}\left\|\dt\phil\right\|_{L^2(0,T;H)}\,+\,\left\|\Bs(\phil)\right\|_{L^\infty(0,T;H)}
\,+\,\bigl\|\Fl(\phil)+C_0\bigr\|_{L^\infty(0,T;\Luno)}\nonumber\\[1mm]
&+\,\left\|\Sl\right\|_{H^1(0,T;V_C^{-\tau})\cap L^\infty(0,T;H)\cap L^2(0,T;\VCt)}
\,+\,\left\|P^{1/2}(\phil)(\Sl-\mul)\right\|_{L^2(0,T;H)}\, \le\, \pier{\hat C_1}\,,
\end{align}
where the constant $\pier{\hat C_1}>0$ is independent of $\alpha\,,\beta\,,\lambda$
\pier{and has the same structure as the \rhs\ of \eqref{bounds1},}  and where $C_0>0$ is a constant such that
$F^\lambda(s)+C_0\ge 0$ for all $s\in\erre$. 
Owing to \gianni{\eqref{yosida4}}, we may take $C_0=\hat c_2$ in our case. 
A~fortiori, \gianni{\eqref{yosida4}} and \eqref{bounds1l} imply that, by choosing a possibly larger
$\pier{\hat C_1}$, we may assume that  
\begin{equation}
\label{bounds2l}
\|\phil\|_{L^\infty(0,T;\VBs)}\,\le\,\pier{\hat C_1}\,.
\end{equation}
Since, by the global Lipschitz continuity of $f_2$, the nonlinearity $F_2$ grows at most quadratically, we
then can also infer the bounds
\begin{equation}
\label{bounds3l}
\|F_1^\lambda(\phil)\|_{L^\infty(0,T;\Luno)}\,+\,\|F_2(\phil)\|_{L^\infty(0,T;\Luno)}\,\le\,\pier{\hat C_1}\,.
\end{equation}
The existence result and the global bound
\eqref{bounds1} then follow from a passage to the limit as $\lambda\searrow0$ in the system \eqref{primal}--\eqref{cauchyl}
and in \eqref{bounds1l}.
\Erem

\vspace{3mm}
In the general case, the equation for $\phi$ is just the variational inequality~\eqref{seconda},
and we cannot write anything that is similar to \eqref{Iseconda},
since no estimate for $f_1(\phi)$ is available.
However, if one reinforces the assumptions on the structure, then
one can recover~\eqref{Iseconda} at least as a differential inclusion.
The crucial condition is the following:
\Bsist
  && \psi(v) \in H
  \aand
  \bigl( B^{2\sigma} v , \psi(v) \bigr) \geq 0,
  \quad \hbox{for every $v\in\VB{2\sigma}$}
  \,\hbox{ and every monotone} \non\\
  &&\mbox{\gianni{and \Lip\ continuous function \,$\psi:\erre\to\erre$ vanishing at the origin}}.
   \label{hpregularity}
\Esist
We notice that this assumption is fulfilled if $B^{2\sigma}=-\Delta$ 
with zero Neumann boundary conditions.
Indeed, in this case it results that $\VB{2\sigma}=\{v\in\Hdue:\ \partial_\pn v=0\}$
and, for every $\psi$ as in \eqref{hpregularity} and $v\in\VB{2\sigma}$, 
we have that $\psi(v)\in\Huno$ (since $v\in\Huno$)~and
\Beq
  \bigl( B^{2\sigma} v, \psi(v) \bigr)
  = \iO (-\Delta v) \, \psi(v)
  = \iO \nabla v \cdot \nabla\psi(v)
  = \iO \psi'(v) |\nabla v|^2
  \geq 0 .
  \non
\Eeq
More generally, in place of the Laplace operator, we can take the principal part
of an elliptic operator in divergence form with \Lip\ continuous coefficients,
provided that the normal derivative is replaced by the conormal derivative.
In any case, we can take the Dirichlet boundary conditions instead of the Neumann boundary conditions,
\gianni{since the functions $\psi$ for which \eqref{hpregularity} is required satisfy $\psi(0)=0$.}

The following result has been proved in \cite[Thm.~2.6]{CGS23}.
\Bthm
\label{Secondeq}
Let the assumptions {\bf (A1)}--{\bf (A4)} be fulfilled,
and assume that $\alpha>0$ and $\beta>0$. If, in addition, \eqref{hpregularity} is satisfied, then
there exist a  solution $\soluz$ to the problem \eqref{prima}--\eqref{cauchy} and some $\,\xi\,$ such that
\begin{align}
&\phi\in L^2(0,T;V_B^{2\sigma})\quad\mbox{and}\quad \xi\in L^2(0,T;H),\\
&\beta\,\dt\phi +B^{2\sigma}\phi+\xi+f_2(\phi)=\mu
\quad\mbox{and}\quad \xi\in f_1(\phi)\,\,\mbox{ a.e. in }\,Q.
\end{align}
Moreover, also $\,\xi\,$ is unique if \eqref{embedAC} holds true, and if we
also assume that
the condition \eqref{ini2} is valid, then the unique solution $\soluz$ and 
the associated $\,\xi\,$ satisfy \eqref{adregmu}--\eqref{adregS} as well as 
\begin{equation}
\label{fullreg}
\phi\in L^\infty(0,T;V_B^{2\sigma})\quad\mbox{and}\quad \xi\in L^\infty(0,T;H).
\end{equation}
\Ethm

\vspace{3mm}
We conclude our preparations with a technical lemma that relates to each other the solutions to 
\eqref{prima}--\eqref{cauchy} for different pairs $(\alpha_i,\beta_i)$, $i=1,2$. 

\Blem
\label{Juergen}
Suppose that {\bf (A1)}--{\bf (A4)} are fulfilled, and let $(\mu_{\alpha_i,\beta_i},\phi_{\alpha_i,\beta_i},
S_{\alpha_i,\beta_i})$ be solutions to \eqref{prima}--\eqref{cauchy} in the sense of Theorem~\gianni{\ref{Wellposed}} 
for the parameters  $(\alpha_i,\beta_i)\in (0,1]$, $i=1,2$. 
Then there is some $\hat M>0$, which only depends on
the global constant\pier{%
$$\widehat K_1 \,\bigl(
    \alpha^{1/2} \norma\muz
    \,+\, \norma{B^\sigma\phiz}
    \,+\, \norma{F(\phiz)}_{\Luno}
    \,+\, \norma\Sz
    + 1  \bigr)
$$
in the \rhs\ of} \eqref{bounds1}, such that, for every $t\in (0,T)$ and every $\delta>0$, we have  
\begin{align}
\label{sirtoby}
&(1-\alpha_1 L-\delta)\int_{Q_t}\bigl|\phi_{\alpha_1,\beta_1}-\phi_{\alpha_2,\beta_2}\bigr|^2\,\le\,
\frac 1{4\delta}\int_{Q_t}\bigl|(\alpha_1\mu_{\alpha_1,\beta_1}+\phi_{\alpha_1,\beta_1})-
(\alpha_2\mu_{\alpha_2,\beta_2}+\phi_{\alpha_2,\beta_2})\bigr|^2 \non\\
&\qquad+\,\hat M\bigl|\alpha_1-\alpha_2\bigr|\,+\,\alpha_1\bigl|\beta_1-\beta_2\bigr|\int_0^t
\|\dt\phi_{\alpha_2,\beta_2}(s)\|\,\|\phi_{\alpha_1,\beta_1}(s)-\phi_{\alpha_2,\beta_2}(s)\|\,ds\,.                            \end{align}  
\Elem  
\Bdim
For convenience, we \pier{set} $\phi_i:=\phi_{\alpha_i,\beta_i}$, $\mu_i:=\mu_{\alpha_i,\beta_i}$, for $i=1,2$. Then we 
multiply \eqref{seconda}, written for $\beta_1,\mu_1,\phi_1$, by $\alpha_1$, insert $v=\phi_2(t)$, 
and add the term $(\phi_1(t),\phi_1(t)-\phi_2(t))$ to both sides of the resulting inequality. We then obtain, almost
everywhere in $(0,T)$, the inequality
\begin{align*}
&(\alpha_1\beta_1\,\dt\phi_1,\phi_1-\phi_2)\,+\,\alpha_1(\Bs\phi_1,\Bs(\phi_1-\phi_2))\,+\,(\phi_1,\phi_1-\phi_2)\\[0.5mm]
&\le\,(\alpha_1\,\mu_1+\phi_1,\phi_1-\phi_2)\,-\,(\alpha_1\,f_2(\phi_1),\phi_1-\phi_2)\,+\,\alpha_1 \pier{\iO} (F_1(\phi_2)-F_1(\phi_1)).
\end{align*}
Similarly, arguing on the inequality for $\beta_2,\mu_2,\phi_2$, we get
\begin{align*}
&(\alpha_2\beta_2\,\dt\phi_2,\phi_2-\phi_1)\,+\,\alpha_2(\Bs\phi_2,\Bs(\phi_2-\phi_1))\,+\,(\phi_2,\phi_2-\phi_1)\\[0.5mm]
&\le\,(\alpha_2\,\mu_2+\phi_2,\phi_2-\phi_1)\,-\,(\alpha_2\,f_2(\phi_2),\phi_2-\phi_1)\,+\,\alpha_2\pier{\iO}(F_1(\phi_1)-F_1(\phi_2)).
\end{align*}
Adding the two inequalities, and rearranging terms, we find that 
almost everywhere in $(0,T)$ it holds the inequality
\begin{align}
\label{yes1}
&\bigl(\alpha_1\beta_1\,\dt\phi_1-\alpha_2\beta_2\,\dt\phi_2,\phi_1-\phi_2)\,+\,\|\phi_1-\phi_2\|^2\,+\,\alpha_1
\bigl\|\Bs(\phi_1-\phi_2)\bigr\|^2\non\\[0.5mm]
&\le \,\bigl((\alpha_1\mu_1+\phi_1)-(\alpha_2\mu_2+\phi_2),\phi_1-\phi_2\bigr)\,-\,(\alpha_1-\alpha_2)(\Bs\phi_2,\Bs(\phi_1-\phi_2))
\non\\[0.5mm]
&\quad\,-\,\alpha_1(f_2(\phi_1)-f_2(\phi_2),\phi_1-\phi_2)\,-\,(\alpha_1-\alpha_2)(f_2(\phi_2),\phi_1-\phi_2)\non\\[0.5mm]
&\quad\,-\,(\alpha_1-\alpha_2)\pier{\iO}(F_1(\phi_1)-F_1(\phi_2))\,.
\end{align}
Now, \pier{recalling} \eqref{hpFdue}, we see that
\begin{equation}
\label{yes2}
-\,\alpha_1(f_2(\phi_1)-f_2(\phi_2),\phi_1-\phi_2)\,\le\,\alpha_1 L\|\phi_1-\phi_2\|^2\,.
\end{equation}
Moreover, we have the identity
\begin{align}
\label{yes3}
&\bigl(\alpha_1\beta_1\,\dt\phi_1-\alpha_2\beta_2\,\dt\phi_2,\phi_1-\phi_2)\non\\[0.5mm]
&=\,\frac{\alpha_1\beta_1}2\,\frac d{dt}\,\|\phi_1-\phi_2\|^2\,+\,\bigl((\alpha_1-\alpha_2)\beta_2\,+\,\alpha_1
(\beta_1-\beta_2)\bigr)(\dt\phi_2,\phi_1-\phi_2)\,.
\end{align}
At this point, we integrate the inequality \eqref{yes1} over $(0,t)$. Omitting two nonnegative terms on the \lhs, invoking 
\eqref{yes2} and \eqref{yes3}, and applying the Cauchy--Schwarz and Young inequalities, we find that
\begin{align}
\label{yes4}
&(1-\alpha_1 L)\int_{Q_t}|\phi_1-\phi_2|^2\,\le\,\delta \int_{Q_t}|\phi_1-\phi_2|^2\,+\,\frac{1}{4\delta}\int_{Q_t}|(\alpha_1\mu_1+\phi_1)-(\alpha_2\mu_2+\phi_2)|^2
\non\\
&+\,|\alpha_1-\alpha_2|\Big(\int_0^t\|\Bs\phi_2(s)\|\,\|B^\sigma(\phi_1(s)-\phi_2(s))\|\,ds\,+\,\int_{Q_t}
\bigl(F_1(\phi_1)+F_1(\phi_2)\bigr)\non\\
&\hspace*{25mm}+\,\int_0^t\bigl(\|f_2(\phi_2(s)\pier{)}\|\,+\,\beta_2\|\dt\phi_2(s)\|\bigr)\,\|\phi_1(s)-\phi_2(s)\|\,ds\Big)\non\\
&+\,\alpha_1\,|\beta_1-\beta_2|\int_0^t\|\dt\phi_2(s)\|\,\|\phi_1(s)-\phi_2(s)\|\,ds\,.
\end{align}
Finally, observe that the expression in the bracket multiplying $\,|\alpha_1-\alpha_2|\,$ is, owing to \eqref{bounds1},
bounded in terms of the constant $\hat K_1$. From this, the assertion follows.
\Edim

\section{The case $\mathbf{\alpha\searrow0,\,\beta>0.}$}
\setcounter{equation}{0}

In order to indicate their dependence on the parameters $\alpha,\,\beta$, we denote in the following 
solution triples of the problem \eqref{prima}--\eqref{cauchy} by $(\muab,\phiab,\Sab)$, for
$\,\alpha,\beta\in [0,1)$.  
In this section, we investigate their asymptotic behavior as $\alpha\searrow0$ and $\beta>0$. 
Obviously, the main difficulty in the limit processes is to pass through the limit in the
nonlinearities, which requires a strong convergence of the arguments $\phiab$, in particular. 
Denoting in the following by $\mathbf{1}$ both the functions that are identically equal to unity
on $\Omega$ or $Q$, we assume, in addition
to the general assumptions {\bf (A1)}--{\bf (A4)}:

\vspace{2mm}\noindent
{\bf (A5)} \,\,At least one of the following three conditions is satisfied:\\[0.5mm]
\hbox to 2em{(i)\hfil} $\lambda_1$\, is positive.\\[0.5mm]
\hbox to 2em{(ii)\hfil} $P(s)\ge P_0$ for all $s\in\erre$ and some fixed $P_0>0$.\\[0.5mm]
\hbox to 2em{(iii)\hfil} \gianni{$\lambda_1=0$ is a simple eigenvalue of~$A$ and ${\bf1}$ is an eigenfunction belonging to $\VBs$;}\\%
\hbox to 2em{\hfil} \gianni{moreover, $D(F_1)=\erre$, and there are constants $\hat c_3>0$ and $\hat c_4\ge 0$ such that}
\begin{equation}
\label{lulu}
\gianni{|s'|\,\le\,\hat c_3\,F_1(s)+\hat c_4 \quad\mbox{\juergen{whenever} \,$s\in\erre\,$ and $\,s'\in f_1(s)$}.}
\end{equation}

\vspace{2mm}
\Brem
\label{RemA5}
The condition {\bf (A5)},(i)  is satisfied
by the standard \gianni{second-order} elliptic operators with zero Dirichlet boundary conditions
(however, also zero mixed and Robin boundary conditions can be considered, 
with proper definitions of the domains of the operators). 
The case {\bf (A5)},(ii) is, unfortunately,
not too realistic in the practical application to tumor growth models, in 
which, usually, $P$~should also attain the value zero.
Finally, we comment on {\bf (A5)},(iii).
\gianni{The condition $\lambda_1=0$ is satisfied, e.g.,
if $A$ is the Laplace operator $-\Delta$ with zero Neumann boundary conditions.
Furthermore, in this case, the eigenvalue $\lambda_1=0$ is simple, 
and the corresponding eigenfunctions are constants, since $\Omega$ is supposed to be connected.
Furthermore, we have $\,\mathbf{1}\in\VBs\,$ for 
many standard elliptic operators with zero Neumann boundary conditions
(and even with zero Dirichlet boundary conditions if $\sigma$ is small,
for instance, if $B=-\Delta$ with $D(B)=\Hdue\cap\Hunoz$ and $\sigma<1/4$).}
Moreover, the condition
\eqref{lulu} excludes the logarithmic and double obstacle potentials,
\gianni{but it still allows $f_1$ to be multi-valued, since it does not require that $F_1$ is differentiable}; 
it is, however, satisfied for a wide class of smooth potentials of \pier{polynomial 
(and even first-order exponential)} type such as~$F_{reg}$. 
\Erem

\gianni{%
\Brem
\label{Meanvalue}
Clearly, we have that 
$$
  \norma{A^\rho v}^2
  = \somma j1\infty |\lambda_j^\rho (v,e_j)|^2
  \geq \lambda_1^{2\rho} \somma j1\infty |(v,e_j)|^2
  = \lambda_1^{2\rho} \, \norma v^2
  \quad \hbox{for every $v\in\VAr $}.
$$
Hence, in the case {\bf(A5)}(i) in which $\lambda_1>0$, the function $v\mapsto\norma{A^\rho v}$ 
is a norm on $\VA r$ that is equivalent to~\eqref{defnormaAr}.
On the contrary, in the case {\bf(A5)}(iii), 
we have $\lambda_1=0$ and the above function is just a seminorm on~$\VA\rho$.
However, the assumptions that $\lambda_1=0$ is a simple eigenvalue 
and that the eigenfunctions are constants imply the Poincar\'e type inequality (cf.\ \cite[Eq.~(3.5)]{CGS18})
\begin{equation}
\label{poincare}
\|v\|\,\le\,\hat c\,\|\Ar v\|\quad\mbox{\pier{for} some $\,\hat c>0$, for every \,$v\in\VAr\,$ with }\,\mbox{mean}(v)=0,
\end{equation} 
where $\mbox{mean}(v)$ denotes the mean value of~$v$.
Then, a standard argument based on \eqref{poincare} and the compactness of the embedding $\VAr\subset H$
yield that the mapping
\begin{equation}
\label{eqnorm}
v\mapsto |\mbox{mean}(v)|\,+\,\|\Ar v\|\quad\mbox{for }\,v\in\VAr,
\end{equation}
defines a norm on $\VAr$ which is equivalent to $\|\cdot\|_{\VAr}$.
\Erem
}%

\vspace{2mm}
In the following, we denote by $C_i$, $i\in \enne$, positive constants that may depend on the data of the
system but not on the parameters $\alpha,\beta, \lambda$.
We suppose that $\beta>0$ is fixed and $\{\alpha_n\}$ is any sequence satisfying $\alpha_n\searrow0$.
In view of the global bounds \eqref{bounds1}, we may without loss of generality assume that there
are functions $\zeta,\xi, \mu_{0,\beta}, \phi_{0,\beta},S_{0,\beta}$ such that, as~$n\to\infty$,
\begin{align}
\label{convb1}
&\alpha_n\,\mu_{\alpha_n,\beta}\to 0\quad\mbox{strongly in }\,L^\infty(0,T;H),\\[1mm] 
\label{convb2}
&\dt\bigl(\alpha_n\,\mu_{\alpha_n,\beta}+\phi_{\alpha_n,\beta}\bigr) \,\to\,\zeta
\quad\mbox{weakly in }\,L^2(0,T,V_A^{-\rho}),\\[1mm]
\label{convb3}
&A^\rho\mu_{\alpha_n,\beta}\,\to\,\xi \quad\mbox{weakly in }\,
 L^2(0,T;H),\\[1mm]
\label{convb4}
&\phi_{\alpha_n,\beta}\,\to\,\phi_{0,\beta}\quad\mbox{weakly-star in }\,H^1(0,T;H)\cap
L^\infty(0,T;\VBs),\\[1mm]
\label{convb5}
&S_{\alpha_n,\beta}\,\to\, S_{0,\beta}\quad\mbox{weakly-star in } \,H^1(0,T;V_C^{-\tau})\cap 
L^\infty(0,T;H)\cap L^2(0,T;\VCt).
\end{align}
Obviously, \eqref{convb1}, \eqref{convb2} and \eqref{convb4} imply that $\,\zeta=\dt\phi_{0,\beta}$. 
We now claim that
the condition {\bf (A5)} implies that, at least for a subsequence,
\begin{equation}
\label{convb6}
\mu_{\alpha_n,\beta}\,\to\,\mu_{0,\beta}\quad\mbox{weakly in \,$L^2(0,T;\VAr)$},
\end{equation}
which entails, in particular, that $\xi=A^\rho\mu_{0,\beta}$. 

This follows directly if $\lambda_1>0$: indeed, 
\gianni{as observed in Remark~\ref{Meanvalue}, the mapping
$\,v\mapsto \|A^\rho v\|\,$ defines a norm on $\VAr$ which is equivalent to $\|\cdot\|_{\VAr}$ in this case,} 
and thus the boundedness of $\,\{\|A^\rho\mu_{\alpha_n,\beta}\|_{L^2(0,T;H)}\}_{n\in\enne}\,$  entails that \eqref{convb6}
holds true at least for a subsequence. 

Suppose next that $\lambda_1=0$ and that {\bf (A5)},(ii) is fulfilled. Then we can
test the equation \eqref{primal} in the Moreau--Yosida approximation, written at the time $\,s$,
 by $v=\mul_{\alpha_n,\beta}(s)$ and integrate
over $(0,t)$ where $t\in (0,T]$. We then obtain the inequality
\begin{align}
\label{P1}
& \gianni{\frac {\alpha_n} 2} \, \bigl\|\mul_{\alpha_n,\beta}(t)\bigr\|^2\,+\,
\int_0^t\bigl(\|\Ar\mul_{\alpha_n,\beta}(s)\|^2+P_0\,\|\mul_{\alpha_n,\beta}(s)\|^2\bigr)ds\nonumber\\
&\le\,\gianni{\frac {\alpha_n} 2} \,\|\muz\|^2\,+\,\int_0^t\!\!\iO\bigl(P(\phil_{\alpha_n,\beta})
\,S^\lambda_{\alpha_n,\beta}\,-\,\dt\phil_{\alpha_n,\beta}\bigr)\,\mul_{\alpha_n,\beta}\,.
\end{align}
Invoking the global bounds \eqref{bounds1l} and Young's inequality, we readily see that the \rhs\
is bounded by an expression of the form
\begin{equation*}
\frac {P_0}2 \int_0^t\|\mul_{\alpha_n,\beta}(s)\|^2\,ds \,+\,C_1\,.
 \end{equation*}
 Therefore, \pier{it turns out that}
$\,\bigl\|\mul_{\alpha_n,\beta}\bigr\|_{L^2(0,T;\VAr)}\,\le\,C_2$. Letting $\lambda\searrow0$, and
invoking the semicontinuity of norms, we then conclude that
\begin{equation}
\label{convb7}
\bigl\|\mu_{\alpha_n,\beta}\bigr\|_{L^2(0,T;\VAr)}\,\le\,C_2 \quad\forall\,n\in\enne,
\end{equation}
which yields  the validity of \eqref{convb6} on a subsequence also in this case.

It remains to show \eqref{convb6} if $\lambda_1=0$ and {\bf (A5)},(iii) is satisfied.
\gianni{We recall that \eqref{yosida3} yields
$\,f_1^\lambda(s)\in f_1(J^\lambda(s))$ for all $s\in\erre$
so that we can \pier{apply} \eqref{lulu} with $s$ replaced by $J^\lambda(s)$ and $s'=f_1^\lambda(s)$. 
Hence, by also using~\eqref{yosida2}}, 
we find for every $\lambda\in\erre$ and $s\in\erre$ the chain of inequalities
\begin{align}
\label{bravopier}
&\gianni{|f_1^\lambda(s)|\,
\le\,\hat c_3 \,F_1(J^\lambda(s))\,+\,\hat c_4\,
\le\,\hat c_3\,F_1^\lambda(s)\,+\,\hat c_4\,.}
\end{align}
Now recall that $f_2$ is globally Lipschitz continuous on $\erre$, whence it follows that $f_2$ grows at most
linearly and $F_2$ grows at most quadratically. Hence, invoking also \eqref{yosida4}, we can infer that, for every $s\in\erre$,
\begin{align}
\label{convb8}
|f^\lambda(s)| \,&\le\,|f_1^\lambda(s)|\,+\,|f_2(s)|\,\le \hat c_3\,F_1^\lambda(s)+\hat c_4+|f_2(s)|\nonumber\\[0.5mm]
&\le\,\hat c_3\,F^\lambda(s)\,+\,\hat c_3\,|F_2(s)|\,+\,\hat c_4+|f_2(s)|\,\le\,\hat c_3\,F^\lambda(s)\,+\,C_3(1+s^2)
\nonumber\\[0.5mm]
&\le\,\bigl(\hat c_3\,+\,C_3\,{\hat c_1}^{\,\,-1}\bigr)\,F^\lambda(s)\,+\,C_4\,.
\end{align}
Therefore, we can conclude from \eqref{bounds1l}
and \eqref{yosida4} the bounds
\begin{align}
\label{convb9}
\bigl\|f^\lambda(\phil_{\alpha_n,\beta})\bigr\|_{L^\infty(0,T;\Luno)}\,+\,
\bigl\|\phil_{\alpha_n,\beta}\bigr\|_{L^\infty(0,T;H)}\,\le\,C_5\,.
\end{align} 
At this point, we insert $v=\pm\mathbf{1}\in\VBs$ in
\eqref{eqsecondal} to find the estimate
\gianni{%
\begin{align*}
\pm\iO\mul_{\alpha_n,\beta}(t)\,
&\le\,C_6\Bigl(\pcol{\beta}\,\|\dt\phil_{\alpha_n,\beta}(t)\|\,
+\,\|\Bs \phil_{\alpha_n,\beta}(t)\|\,\|\Bs\mathbf{1}\|\,
+\,\bigl\|f^\lambda(\phil_{\alpha_n,\beta})(t)\bigr\|_{\Luno}\Bigr)\\[0.5mm]
&\le\,C_6\Bigl(\pier{\beta}\,\|\dt\phil_{\alpha_n,\beta}(t)\|\,
+\,\|\phil_{\alpha_n,\beta}(t)\|_{\VBs}\,\|\Bs\mathbf{1}\|\,
+\,C_5\Bigr)\,,
\end{align*}
}%
which, owing to \eqref{bounds1l} and \gianni{\eqref{convb9}}, then shows that
\begin{equation*}
\bigl\|\mbox{mean} (\mul_{\alpha_n,\beta})\bigr\|_{L^2(0,T)}\,\le\,C_7\,\pier{(1+\beta)}\,.
\end{equation*}
Combining this with \eqref{bounds1l}, and recalling the 
equivalence of the norms \eqref{eqnorm} and \linebreak $\norma\cpto_{\VAr}$ \gianni{given in Remark~\ref{Meanvalue}}, 
we have finally shown that the sequence $\,\{\|\mul_
{\alpha_n,\beta}\|_{L^2(0,T;\VAr)}\}_{n\in\enne}\,$ is bounded. 
Passage to the limit as
$\lambda\searrow0$, and the semicontinuity of norms,  then yield that also 
$\,\{\|\mu_{\alpha_n,\beta}\|_{L^2(0,T;\VAr)}\,\}_{n\in\enne}\,$ is bounded. 
With this, 
we can conclude the validity of \eqref{convb6} on a subsequence also in this case.

\vspace{2mm}
With \eqref{convb6} shown for all of the cases considered in {\bf (A5)}, we can continue our analysis.  
At first, \pier{thanks to} \eqref{convb4}, \eqref{convb5}, and 
\pier{known compactness results (see, e.g., \cite[Sect.~8, Cor.~4]{Simon}),}
we may without loss of generality assume that
\begin{align}
\label{convb10}
\phi_{\alpha_n,\beta}\,&\to\,\phi_{0,\beta}\quad\mbox{\pier{strongly in \,$C^0([0,T];H)$}},\\[0.5mm]
\label{convb11}
S_{\alpha_n,\beta}\,&\to\,S_{0,\beta}\quad\mbox{strongly in \,$L^2(0,T;H)$.}
\end{align}
Then, by the Lipschitz continuity of both $P$ and $f_2$, 
\begin{align}
\label{convb12}
P(\phi_{\alpha_n,\beta})\,&\to\,P(\phi_{0,\beta})\quad\mbox{\pier{strongly in \,$C^0([0,T];H)$}},\\[0.5mm]
\label{convb13}
f_2(\phi_{\alpha_n,\beta})\,&\to\,f_2(\phi_{0,\beta})\quad\mbox{\pier{strongly in \,$C^0([0,T];H)$}}.
\end{align}
Next, we observe that the convergence properties \eqref{convb6}, \eqref{convb11}, and \eqref{convb12} imply that
\gianni{%
\Beq
  P(\phi_{\alpha_n,\beta})\bigl(S_{\alpha_n,\beta}-\mu_{\alpha_n,\beta}\bigr)
  \to
  P(\phi_{0,\beta})\bigl(S_{0,\beta}-\mu_{0,\beta}\bigr)
  \quad \hbox{weakly in $\LQ1$\,.}
  \non
\Eeq
On the other hand, $P(\phi_{\alpha_n,\beta})\bigl(S_{\alpha_n,\beta}-\mu_{\alpha_n,\beta}\bigr)$
is bounded in $\L2H$ due to \eqref{bounds1l}, since $P$ is bounded.
Hence, we deduce that
\Beq
  P(\phi_{\alpha_n,\beta})\bigl(S_{\alpha_n,\beta}-\mu_{\alpha_n,\beta}\bigr)
  \to
  P(\phi_{0,\beta})\bigl(S_{0,\beta}-\mu_{0,\beta}\bigr)
  \quad \hbox{weakly in $\L2H$\,.}
  \label{convb14}
\Eeq
}%

Now we are in a position to take the limit as $n\to\infty$ in the time-integrated versions of \eqref{prima}
and \eqref{terza}, respectively, written with time-dependent test functions $v$. 
We then obtain that the triple
$(\mu,\phi,S):=(\mu_{0,\beta},\phi_{0,\beta},S_{0,\beta})$ satisfies \eqref{prima} for $\alpha=0$, \eqref{terza}, and 
the initial conditions \gianni{\eqref{cauchyazero}}.

It remains to show the validity of \eqref{seconda} or of its time-integrated version \eqref{intseconda}.
To this end, notice \pier{that the convex functional $v\mapsto \iO F_1(v) $, extended with value $+\infty$ whenever $F_1(v) \not\in \Luno$, is proper, convex and lower semicontinuous in $H$. 
Hence, the convergence \eqref{convb10} and the bound \eqref{bounds1} imply that}
\begin{align}
\non 
0\,\le\,\pier{\iO F_1(\phi_{0,\beta}(t) )\,\le \,\liminf_{n\to\infty} \iO F_1(\phi_{\alpha_n,\beta}(t))\leq C_8} \quad \pier{\mbox{for every }\, t\in [0,T]},
\end{align}
\pier{for some uniform constant $C_8$. 
It therefore follows that $\,F_1(\phi_{0,\beta})\in \L\infty{\Luno} $\pcol{,} and Fatou's lemma allows us to infer that}
\begin{equation}
\label{convb16}
0\,\le\,\int_Q F_1(\phi_{0,\beta})\,\le\,\liminf_{n\to\infty}\int_Q F_1(\phi_{\alpha_n,\beta})\,<\,+\infty.
\end{equation} 
Moreover, the quadratic form 
$\,v\mapsto \int_0^T(\Bs v(t),\Bs v(t))\,dt\,$ 
is weakly sequentially lower \gianni{semicontinuous} on $L^2(0,T;\VBs)$, which entails that
\begin{equation}
\label{convb17}
\int_0^T \bigl(\Bs \phi_{0,\beta}(t),\Bs \phi_{0,\beta}(t)\bigr)\,dt\,\le\,\liminf_{n\to\infty}
\int_0^T \bigl(\Bs \phi_{\alpha_n,\beta}(t),\Bs \phi_{\alpha_n,\beta}(t)\bigr)\,dt\,.
\end{equation}
Using all of the above convergence results, we can therefore conclude that, for every $v\in L^2(0,T;\VBs)$,
\begin{align}
\label{convb18}
&\intQ F_1(\phi_{0,\beta})\,+\,\int_0^T\bigl(\Bs \phi_{0,\beta}(t),\Bs 
(\phi_{0,\beta}(t)-v(t)\bigr)\,dt\nonumber\\
&\le\,\liminf_{n\to\infty}\,\biggl(\intQ F_1(\phi_{\alpha_n,\beta})\,+
\int_0^T\bigl(\Bs \phi_{\alpha_n,\beta}(t),\Bs (\phi_{\alpha_n,\beta}(t)-v(t))\bigr)\,dt\biggr)\nonumber\\
&\le\,\liminf_{n\to\infty}\,\biggl(\int_0^T \bigl(\mu_{\alpha_n,\beta}(t)-f_2(\phi_{\alpha_n,\beta}(t))
-\beta\,\dt\phi_{\alpha_n,\beta}(t),\phi_{\alpha_n,\beta}(t)-v(t)\bigr)\,dt\, \pier{{}+\,\intQ F_1(v)}\biggr)\nonumber\\
&=\,\int_0^T\bigl(\mu_{0,\beta}(t)-f_2(\phi_{0,\beta}(t))-\beta\,\dt\phi_{0,\beta}(t),\phi_{0,\beta}(t)-v(t)\bigr)\,dt\, \pier{{}+\,\intQ F_1(v),}
\end{align} 
which shows the validity of \eqref{intseconda} for $(\mu,\phi,S)=(\mu_{0,\beta},\phi_{0,\beta},S_{0,\beta})$. From the above
analysis, we can conclude the following existence and convergence result.

\Bthm
\label{Atozero}
Suppose that the conditions {\bf (A1)}--{\bf (A5)} are fulfilled, let $\beta>0$ be fixed and $\{\alpha_n\}_{n\in\enne}\subset (0,1]$ be a sequence such that $\alpha_n\searrow 0$. Then there are a subsequence 
$\{\alpha_{n_k}\}_{k\in\enne}$
and functions $(\mu_{\alpha_{n_k},\beta}, \phi_{\alpha_{n_k},\beta},S_{\alpha_{n_k},\beta})$, which solve the system 
\eqref{prima}--\eqref{cauchy} for $\alpha=\alpha_{n_k}$ in the sense of Theorem~\gianni{\ref{Wellposed}}, such that 
there is a triple $(\mu_{0,\beta},\phi_{0,\beta},S_{0,\beta})$ with the following properties:
\begin{align}
\label{condtmuphi}
&\dt\bigl(\alpha_{n_k}\,\mu_{\alpha_{n_k},\beta}\,+\,\phi_{\alpha_{n_k},\beta}\bigr)\,\to\,
\dt\phi_{0,\beta} \quad\mbox{weakly in }\,L^2(0,T;V_A^{-\rho}),\\[1mm]
\label{conbmu}
&\mu_{\alpha_{n_k},\beta}\,\to \,\mu_{0,\beta}\quad\mbox{weakly-star in }\,L^2(0,T;\VAr),\\[1mm]
&\phi_{\alpha_{n_k},\beta}\,\to \,\phi_{0,\beta}\quad\mbox{weakly-star in }\,H^1(0,T;H)\cap
 L^\infty(0,T;\VBs),\\[1mm]
&S_{\alpha_{n_k},\beta}\,\to\,S_{0,\beta}\quad\mbox{weakly-star in }\,H^1(0,T;V_C^{-\tau})\cap
L^\infty(0,T;H)\cap L^2(0,T;\VCt).   
\end{align}
In addition, $\,F_1(\phi_{0,\beta})\in \pier{\L\infty{\Luno}}$, and $\,(\mu_{0,\beta},\phi_{0,\beta},S_{0,\beta})$\, 
solves the system \eqref{prima}--\eqref{terza} for $\alpha=0$ 
and satisfies the initial conditions \gianni{\eqref{cauchyazero}}. 
Finally, it holds the additional regularity
\begin{equation}
\label{regmub}
\mu\in L^2(0,T;V_A^{2\rho}).
\end{equation} 
\Ethm

\Bdim
Except for \eqref{regmub}, everything was already proved above. The validity of \eqref{regmub} follows directly from 
comparison in \pier{\eqref{prima}}, since, owing to the boundedness of $\,P$, we have  $\,P(\phi_{0,\beta})\bigl(S_{0,\beta}-\mu_{0,\beta}\bigr)-\dt\phi_{0,\beta}
\in L^2(0,T;H)$. 
\Edim 

Next, we give a regularity result that resembles the corresponding results \eqref{adregmu}--\eqref{adregS}
in  Theorem~\gianni{\ref{Wellposed}} for the 
case when both $\alpha>0$ and $\beta>0$. Note that we cannot expect the same regularity  here, since a vanishing $\,\alpha\,$ entails a loss of coercivity with respect to the solution component $\mu$.

\vspace{2mm}
\Bthm
\label{RegAzero}
Suppose that {\bf (A1)}--{\bf (A4)}, \eqref{ini2}, \eqref{embedAC}, and at least one of the two conditions
\begin{equation}
\label{coercmu}
{\rm (i) }\,\,\lambda_1>0, \quad\mbox{and}\quad {\rm (ii) }\,\,P(s)\ge P_0>0\,\mbox{ for all }\,s\in\erre,
\end{equation}
\juergen{are} fulfilled. Then the solution $(\mu_{0,\beta},\phi_{0,\beta},S_{0,\beta})$ 
established in Theorem~\gianni{\ref{Atozero}} enjoys the additional regularity
\begin{align}
\label{addregmub}
&\mu_{0,\beta}\in L^\infty(0,T;V_A^{2\rho}),\\[0.5mm]
\label{addregphib}
&\phi_{0,\beta}\in W^{1,\infty}(0,T;H)\cap H^1(0,T;\VBs),\\[0.5mm]
\label{addregSb}
&S_{0,\beta}\in H^1(0,T;H)\cap L^\infty(0,T;\VCt)\cap L^2(0,T;V_C^{2\tau}).
\end{align} 
\Ethm

\Bdim
Let, for convenience, $(\mu,\phi, S):=(\mu_{0,\beta},\phi_{0,\beta},S_{0,\beta})$.
We only give a formal proof of the assertion based on the Moreau--Yosida approximation,
which is for $\lambda>0$ given by the system \eqref{primal} with $\alpha=0$, \eqref{eqsecondal}
\gianni{(in~place of \eqref{secondal})},
\eqref{terzal}, together with the initial condition \gianni{\eqref{cauchyazero}}. 
For a rigorous proof,
one would have to carry out the following arguments on the level 
of the time-discretized version introduced in~\cite{CGS23}. 
Since this requires a considerable writing effort without bringing
new insights in comparison with the calculations in~\cite{CGS23}, we prefer to argue formally, here. 
To this end, we differentiate \eqref{eqsecondal} with respect to $t$ and
\pier{take $v=\dt\phil(t)$} in the resulting equation. 
In addition, we insert \pier{$v=\dt\mul(t)$} in \eqref{primal}, 
add the two resulting equations, and integrate their sum over $(0,t)$ where $t\in(0,T]$. 
Noting that the two terms involving $\,\dt\phil\,\dt\mul\,$ cancel each other, we arrive at the identity
\begin{align}
\label{susi1}
&\frac \beta 2\,\|\dt\phil(t)\|^2\,+\,\frac 12 \,\|\Ar\mul(t)\|^2\,+\int_{Q_t}|\Bs(\dt\phil)|^2
\,+\int_{Q_t}(f_1^\lambda)'(\phil)\,|\dt\phil|^2\nonumber\\
&=\,\int_{Q_t} P(\phil)\,(\Sl-\mul)\,\dt\mul \,+\,
\frac \beta 2\,\|\dt\phil(0)\|^2\,+\,\frac 12\,\|\Ar \muz\|^2\,-\int_{Q_t} f_2'(\phil)\,|\dt\phil|^2\,,
\end{align}
where, due to the general assumptions, all of the terms on the \lhs\ are nonnegative and the last term
on the \rhs\ is from \eqref{bounds1l} already known to be bounded independently of $\lambda$.
Now observe that, by formal insertion of $\,\dt\phil(0)$\, in \eqref{eqsecondal} for $t=0$, it follows
from \eqref{ini2} that
\begin{align}
\label{susi2}
\beta\,\|\dt\phil(0)\|^2\,&=\,\left(-B^{2\sigma}\phi_0-f^\lambda(\phi_0)+\mu_0,\dt\phil(0)\right)\,
\le\,\frac \beta 2\,\|\dt\phil(0)\|^2\,+\,C_1,
\end{align}
where, here and in the remainder of this proof, we denote by $C_i$, $i\in\enne$, positive constants
that do not depend on $\lambda$. 
Next, an integration by parts yields that
\begin{align}
\label{susi3}
&-\int_{Q_t} P(\phil)\,\mul\,\dt\mul\,=\,\frac 12\iO\bigl(P(\phiz)\,|\muz|^2-
 P(\phil(t))\,|\mul(t)|^2\bigr) \,+\,\frac 12\int_{Q_t} P'(\phil)\dt\phil\,|\mul|^2 \non\\
&\le\,C_2\,-\,\frac 12\iO P(\phil(t))\,|\mul(t)|^2\,+\,C_3\int_0^t \|\mul(s)\|_{L^4(\Omega)}^2
\,\|\dt\phil(s)\|\,ds\non\\
&\le\,\gianni{C_2}\,-\,\frac 12\iO P(\phil(t))\,|\mul(t)|^2\,
+\,\gianni{C_4}\int_0^t \|\mul(s)\|_{\VAr}^2\,\|\dt\phil(s)\|^2\,ds  \,,
\end{align}
where we used \gianni{H\"older's inequality and \eqref{embedAC}.
\pier{Note that, by virtue of \eqref{hpP}, the second term on the \rhs\ of \eqref{susi3} 
is nonpositive so that it can be moved with the right sign on the \lhs\ of \eqref{susi1}.}
Moreover, we notice that $\,\|\mul\|_{L^2(0,T;\VAr)}$ is uniformly bounded with respect to~$\lambda$
as shown in the proof of Theorem~\ref{Atozero},
and we will account for this information in applying Gronwall's lemma.}
Moreover, \gianni{integrating by parts and} using also the already known bounds~\eqref{bounds1l},
\gianni{the inequality $P\leq P^{1/2}(\sup P^{1/2})$ and the Young inequality},
we infer that 
\begin{align}
\label{susi4}
&\int_{Q_t} P(\phil)\,\Sl\,\dt\mul\non\\
&=\,\iO\bigl(P(\phil(t)\,\Sl(t)\,\mul(t)-\pier{P(\phiz)\,\Sz\,\muz}\bigr)
\,-\,\int_{Q_t}\,\bigl(P'(\phil)\,\dt\phil\,\Sl+P(\phil)\,\dt\Sl\bigr)\mul
 \non\\
&\le\,\gianni{C_5}\,
+\,\frac 14\iO P(\phil(t))\,|\mul(t)|^2\,
\pier{+\,\gianni{C_6}\int_0^t \|\dt\phil(s)\|\, \|\Sl(s)\|_{\VCt}\, \|\mul(s)\|_{\VAr}\,ds} 
\non\\
&\quad
+\,\frac 12\int_{Q_t}|\dt\Sl|^2
\pier{{}+C_7 \intQt P(\phil)|\mul|^2}
\,.
\end{align}
\pier{Note that the third term on the \rhs\ can be treated for instance as 
\begin{align}
\label{susipier}
&
\gianni{C_6}\int_0^t \|\dt\phil(s)\|\, \|\Sl(s)\|_{\VCt}\, \|\mul(s)\|_{\VAr}\,ds
\non\\
&\le
\,\gianni{C_6}\int_0^t  \|\Sl(s)\|_{\VCt}\, \|\mul(s)\|_{\VAr}\,\tonde{1+ \|\dt\phil(s)\|^2}\,
ds\pcol{,}
\end{align}%
and both $\,\|\Sl\|_{L^2(0,T;\VCt)}$ and $\,\|\mul\|_{L^2(0,T;\VAr)}$ are uniformly bounded with respect to~$\lambda$ (cf.~\eqref{bounds1})}.
Finally, we test \eqref{terzal} by $\dt\Sl(t)$ and integrate over $(0,t)$. Then we obtain 
\begin{align}
\label{susi5}
&\int_{Q_t}|\dt\Sl|^2\,+\,\frac 12\,\|\Ct\Sl(t)\|^2\,\le\,\frac 12\,\|\Ct S_0\|^2
\,+\,\int_0^t\bigl(z^\lambda(s),\dt\Sl(s)\bigr)\,ds\non\\
&\le \,\juergen{\frac 12\,\|\Ct S_0\|^2\,+\,}\frac 14\int_{Q_t}|\dt\Sl|^2\,+\,\int_{Q_t}|z^\lambda|^2,
\end{align}
where $\,z^\lambda:=P(\phil)(\mul-\Sl)\,$ is already known 
\gianni{to be bounded in $L^2(0,T;H)$, 
independently of~$\lambda$, by \eqref{bounds1l} and the boundedness of~$P$}. 
Combining \eqref{susi1}--\eqref{susi5}, and invoking Gronwall's lemma, we have therefore
shown the estimate
\begin{align}
\label{susi6}
&\|\dt\phil\|_{L^\infty(0,T;H)\cap L^2(0,T;\VBs)}^2\,+\,
\|\Sl\|_{H^1(0,T;H)\cap L^\infty(0,T;\VCt)}^2\non\\[1mm]
&+\,\sup_{t\in (0,T)} \Big( \|\Ar\mul(t)\|^2\,+\,\iO P(\phil(t))|\mul(t)|^2\Big)\,\le\,\pier{C_8}.
\end{align}
In particular, it follows from \eqref{coercmu} \gianni{(see Remark~\ref{Meanvalue})} that 
\begin{equation}
\label{susi7}
\|\mul\|_{L^\infty(0,T;\VAr)}\,\le\,\pier{C_9}.
\end{equation}
It remains to show the boundedness of $\,\Sl\,$ in $L^2(0,T;V_C^{2\tau})$ and 
of $\,\mul\,$ in $L^\infty(0,T;V_A^{2\rho})$. But this follows
immediately from \eqref{terzal} and \eqref{primal}, respectively,
 by comparison. At this point, we take the limit as $\lambda\searrow0$ 
and invoke the semicontinuity of norms to \pier{infer} that the derived bounds
are valid also in the limit. This concludes the proof of the assertion.
\Edim
   
\vspace{3mm}
\Brem
It is also possible to prove a uniqueness result for the case $\alpha=0$, $\beta>0$,
under restrictive additional assumptions. Since the related analysis requires a major detour in the line
of argumentation and \pier{is carried out in detail in the recent} paper~\cite{CGS25},
we do not present it here. Note also that in the case $P\equiv 0$ the system \eqref{prima}, \eqref{seconda} coincides for
$\alpha=0$ and $\beta\ge 0$ with the system that has recently been studied by the present authors in a series
of papers (see\pier{\cite{CGS18,CGS19,CGS21}}); for precise results in this
much simpler case, in which \eqref{prima}, \eqref{seconda} decouple from \eqref{terza}, we refer to these works. 
\Erem

\section{The case $\mathbf{\alpha>0,\,\beta\searrow0.}$}
\setcounter{equation}{0}

In this section, we investigate the asymptotic behavior of the solutions $(\muab,\phiab,\Sab)$
as $\alpha>0$ and $\beta\searrow0$. In this case, an additional coercivity condition for $\mu$
like \eqref{coercmu} is not necessary. Instead, the main difficulty is to establish a 
strong convergence for the phase variable $\phi$. Indeed, we have to make the following additional assumption:

\vspace{2mm}\noindent
{\bf (A6)} \,\,It holds $\,\,\alpha L<1$.

\vspace{2mm}
\noindent Now, let $\{\beta_n\}_{n\in\enne}\subset (0,1]$ be any
sequence such that $\beta_n\searrow0$. Then, according to the global bound \eqref{bounds1}, we may 
without loss of generality assume the existence of functions $\zeta,\mu_{\alpha,0},\phi_{\alpha, 0},
S_{\alpha, 0}$ such that, \pier{at least for a subsequence} as $n\to\infty$,
\begin{align}
\label{conva1}
&\beta_n\,\dt\phi_{\alpha,\beta_n}\to 0\quad\mbox{strongly in }\,L^2(0,T;H),\\[1mm] 
\label{conva2}
&\dt\bigl(\alpha\,\mu_{\alpha,\beta_n}+\phi_{\alpha,\beta_n}\bigr) \,\to\,\zeta
\quad\mbox{weakly in }\,L^2(0,T,V_A^{-\rho}),\\[1mm]
\label{conva3}
&\mu_{\alpha,\beta_n}\,\to\,\mu_{\alpha, 0}\quad\mbox{weakly-star in }\,
 L^\infty(0,T;H)\cap L^2(0,T;\VAr),\\[1mm]
\label{conva4}
&\phi_{\alpha,\beta_n}\,\to\,\phi_{\alpha,0}\quad\mbox{weakly-star in }\,L^\infty(0,T;\VBs),\\[1mm]
\label{conva5}
&S_{\alpha,\beta_n}\,\to\, S_{\alpha,0}\quad\mbox{weakly-star in } \,H^1(0,T;V_C^{-\tau})\cap L^\infty(0,T;H)\cap
L^2(0,T;\VCt).
\end{align}
Now, combining \eqref{conva2}--\eqref{conva4}, we see that 
$\,\zeta=\dt\left(\alpha\,\mu_{\alpha,0}+\phi_{\alpha,0}\right)$,
and we infer from \eqref{conva3} and \eqref{conva4} that $\,\{\alpha\,\mu_{\alpha,\beta_n}+\phi_{\alpha,\beta_n}\}_{n\in\enne}\,$
is also bounded in the Banach space $\,L^2(0,T;\VAr+\VBs)$, where 
\begin{align}
\label{XpY1}
&\VAr+\VBs\,:=\,\{v+w:\,v\in\VAr\,\,\mbox{ and }\,\,w\in\VBs\},\quad\mbox{and}\\[1mm]
\label{XpY2}
&\|y\|_{\VAr+\VBs}\,:\,=\inf\,\{\|v\|_{\VAr}\,+\,\|w\|_{\VBs}\,:\,v\in \VAr, \,\,\,
w\in\VBs,\,\,\mbox{ and }\,\,
y=v+w\}.
\end{align}
Since both $\VAr$ and $\VBs$ are compactly embedded in $H$, 
so is $\VAr+\VBs$, and we can infer from the Aubin--Lions compactness lemma 
\pier{(see, e.g., \cite[Thm.~5.1, p.~58]{Lions})} that
\begin{equation}
\label{conva6}
\alpha\,\mu_{\alpha,\beta_n}+\phi_{\alpha,\beta_n}\,\to\,\alpha\,\mu_{\alpha,0}+\phi_{\alpha,0}
\quad\mbox{strongly in }\,L^2(0,T;H).
\end{equation}
 Next, we aim at showing that $\{\phi_{\alpha,\beta_n}\}_{n\in\enne}$ is a Cauchy sequence in $L^2(0,T;H)$,
which would imply that, \pier{possibly taking another subsequence}, 
\begin{align}
\label{conva7}
&\phi_{\alpha,\beta_n}\,\to\,\phi_{\alpha,0}\quad\mbox{strongly in \,$L^2(0,T;H)\,$ and}
\non\\
& \quad \pier{\phi_{\alpha,\beta_n}(t) \,\to\,\phi_{\alpha,0}(t) \quad\mbox{strongly in \,$H$,
 for a.e.\ $t\in (0,T)$.}}
\end{align}
To prove the claim, we employ Lemma~\gianni{\ref{Juergen}} with the special choice $\,(\alpha,\beta_n)\,$ and
$(\alpha,\beta_m)$, where $\,n>m\,$ so that $0<\beta_n <\beta_m$.
Thanks to \eqref{sirtoby}, we have, for every $\delta>0$,
\begin{align}
\label{sirtobya}
&(1-\alpha L-\delta)\int_{\pier{Q}}\bigl|\phi_{\alpha,\beta_n}-\phi_{\alpha,\beta_m}\bigr|^2\,\le\,
\frac 1{4\delta}\int_{\pier{Q}}\bigl|(\alpha\,\mu_{\alpha,\beta_n}+\phi_{\alpha,\beta_n})-
(\alpha\,\mu_{\alpha,\beta_m}+\phi_{\alpha,\beta_m})\bigr|^2 \non\\
&\qquad+\,\alpha\bigl|\beta_n-\beta_m\bigr|\int_0^{\pier{T}}
\|\dt\phi_{\alpha,\beta_m}(s)\|\,\|\phi_{\alpha,\beta_n}(s)-\phi_{\alpha,\beta_m}(s)\|\,ds\,.                 
\end{align}
Now observe that $\,|\beta_n-\beta_m|\,=\,(\beta_m-\beta_n)\,\le\,\beta_m\,$. 
Moreover, $\juergen{\left\{\phi_{\alpha,\beta_n}-\phi_{\alpha,\beta_m}\right\}}$ is bounded in $\L2H$. Hence,
by virtue of \eqref{conva1} and \eqref{conva6}, the \rhs\ of \eqref{sirtobya} converges to zero as
$n>m$ and $m\to\infty$. Therefore, choosing $\delta\in (0,1-\alpha L)$, we conclude from \eqref{sirtobya} that
the above claim is valid. We thus may assume that \eqref{conva7} holds true. But this implies that also 
\begin{align}
\label{conva8}
&\mu_{\alpha,\beta_n}\,\to\,\mu_{\alpha,0}\quad\mbox{strongly in }\,L^2(0,T;H),\\[2mm]
\label{conva9}
&P(\phi_{\alpha,\beta_n})\,\to\,P(\phi_{\alpha,0}) \quad \hbox{and} \quad f_2(\phi_{\alpha,\beta_n})
\,\to\,f_2(\phi_{\alpha,0}) 
\non\\
&\qquad\mbox{both strongly in \,$L^2(0,T;H)$},
\end{align}
using \gianni{\eqref{conva6} and the Lipschitz continuity of $P$ and~$f_2$}. 
Moreover, the Aubin--Lions lemma yields that also
\begin{align} 
\label{conva10}
&S_{\alpha,\beta_n}\,\to\,S_{\alpha,0}\quad\mbox{strongly in }\,L^2(0,T;H),
\end{align}
and, as in \gianni{\eqref{convb14}}, it is readily verified that 
\begin{equation}\label{conva11}
\gianni{%
  P(\phi_{\alpha,\beta_n})(S_{\alpha,\beta_n}-\mu_{\alpha,\beta_n})
  \to P(\phi_{\alpha,0})(S_{\alpha,0}-\mu_{\alpha,0})
  \quad \hbox{weakly in $\L2H$}}.
\end{equation}

Now we are in a position to take the limit as $n\to\infty$ in the time-integrated versions of \eqref{prima}
and \eqref{terza}, respectively, written with time-dependent test functions. We then obtain that the triple
$(\mu,\phi,S):=(\mu_{\alpha,0},\phi_{\alpha,0},S_{\alpha,0})$ satisfies \eqref{prima} and \eqref{terza}, and 
\eqref{conva2} entails that $\,\alpha\,\mu_{\alpha,\beta_n}+\phi_{\alpha,\beta_n}\,\to
\,\alpha\,\mu_{\alpha,0}+\phi_{\alpha,0}$ weakly in $\,C^0([0,T];V_A^{-\rho})\,$, which shows, in particular,
that $\,(\alpha\,\mu_{\alpha,0}+\phi_{\alpha,0})(0)=\alpha\,\muz+\phiz$, \gianni{i.e., the first of~\eqref{cauchyaposbzero}}. 
At the same time, we conclude from 
\eqref{conva5} the weak convergence $\,S_{\alpha,\beta_n}\,\to\,S_{\alpha,0}\,$ in $\,C^0([0,T];V_C^{-\tau})$;
\gianni{therefore, we also have the second of the initial conditions~\eqref{cauchyaposbzero}}.
It remains to show the validity of \eqref{seconda} or its time-integrated 
version \eqref{intseconda}, for $\beta=0$.
To this end, notice that \pier{\eqref{bounds1}, \eqref{conva7}  and the lower semicontinuity of the functional $v\mapsto \iO F_1(v) $ in $H$ imply that
\begin{align}
\label{conva12}
0\,\le\,\pier{\iO F_1(\phi_{\alpha,0}(t) )\,\le \,\liminf_{n\to\infty} \iO F_1(\phi_{\alpha,\beta_n}(t))\leq C} \quad \pier{\mbox{for a.e. }\, t\in (0,T)},
\end{align}
for some constant $C$ independent of $\beta_n$. 
Thus, it follows that $\,F_1(\phi_{\alpha,0})\in \L\infty{\Luno} $ and, by Fatou's lemma,}
\begin{equation}
\label{conva13}
0\,\le\,\int_Q F_1(\phi_{\alpha,0})\,\le\,\liminf_{n\to\infty}\int_Q F_1(\phi_{\alpha,\beta_n})\,<\,+\infty.
\end{equation}
Moreover, the quadratic form 
$\,v\mapsto \int_0^T(\Bs v(t),\Bs v(t))\,dt\,$ is weakly sequentially lower \gianni{semicontinuous} on $L^2(0,T;\VBs)$. 
Therefore, a similar calculation (which needs no repetition here) as in \eqref{convb16} yields the validity of \eqref{intseconda}.

In conclusion, we have the following result.

\Bthm 
\label{Btozero}
Suppose that the conditions {\bf (A1)}--{\bf (A4)} and {\bf (A6)}  are fulfilled.
\gianni{Moreover, let $\alpha>0$ and 
 $\{\beta_n\}\subset (0,1)$ be} a sequence with $\,\beta_n\searrow0$ as $n\to\infty$.
 Then there are a subsequence 
$\{\beta_{n_k}\}_{k\in\enne}$
and functions $(\mu_{\alpha,\beta_{n_k}}, \phi_{\alpha,\beta_{n_k}},S_{\alpha,\beta_{n_k}})$, 
which solve the system 
\eqref{prima}--\eqref{cauchy} for $\beta=\beta_{n_k}$ in the sense of Theorem~\gianni{\ref{Wellposed}}, and a triple 
$(\mu_{\alpha,0},\phi_{\alpha,0},S_{\alpha,0})$ with the following properties:
\begin{align}
\label{coka1}
&\beta_{n_k}\,\dt\phi_{\alpha,\beta_{n_k}}\to 0\quad\mbox{strongly in }\,L^2(0,T;H),\\[1mm] 
\label{coka2}
&\dt\bigl(\alpha\,\mu_{\alpha,\beta_{n_k}}+\phi_{\alpha,\beta_{n_k}}\bigr) \,\to\,\dt\bigl(
\alpha\,\mu_{\alpha,0}+\phi_{\alpha,0}\bigr)
\quad\mbox{weakly in }\,L^2(0,T,V_A^{-\rho}),\\[1mm]
\label{coka3}
&\mu_{\alpha,\beta_{n_k}}\,\to\,\mu_{\alpha, 0}\quad\mbox{weakly-star in }\,
 L^\infty(0,T;H)\cap L^2(0,T;\VAr)\\
 &\hspace*{30mm}\mbox{\pier{and strongly in }}\,L^2(0,T;H),\\[1mm]
\label{coka4}
&\phi_{\alpha,\beta_{n_k}}\,\to\,\phi_{\alpha,0}\quad\mbox{weakly-star in }\,L^\infty(0,T;\VBs)
\,\mbox{ and strongly in }\,L^2(0,T;H),\\[1mm]
\label{coka5}
&S_{\alpha,\beta_{n_k}}\,\to\, S_{\alpha,0}\quad\mbox{weakly-star in } \,H^1(0,T;V_C^{-\tau})\cap L^\infty(0,T;H)\cap
L^2(0,T;\VCt).
\end{align}
In addition, $\,F_1(\phi_{\alpha,0})\in \pier{\L\infty{\Luno}}$, 
and $\,(\mu_{\alpha,0},\phi_{\alpha,0},S_{\alpha,0})$\, solves the system \eqref{prima}--\eqref{terza} for
$\beta=0$ and satisfies the initial conditions \gianni{\eqref{cauchyaposbzero}}. 
\Ethm

\vspace{2mm}
It seems to be difficult to derive additional regularity results for $\alpha>0$ and $\beta=0$,
\gianni{and we give a comment on this in the forthcoming Remark~\ref{RemAposBzero}. 
However, we can \pier{show} a more important uniqueness result.}  
To this end, we need to make a compatibility assumption that strongly relates the operators
$\Ar$ and $\Bs$ to each other. We have the following result.

\Bthm
\label{UniqAzero}
\juergen{Assume, in addition to {\bf (A1)}--{\bf (A4)} and {\bf (A6)}, that the following embeddings are continuous:}
\Beq
  \VBs \subset \VAr \,, \quad
  \VAr \subset \Lx4 \,, \quad
  \VBs \subset \Lx4 
  \aand
  \VCt \subset \Lx4 \,.
  \label{embedABC}
\Eeq
Then the solution to the system \eqref{prima}--\eqref{terza}, \gianni{\eqref{cauchyaposbzero}} for $\alpha>0$ and $\beta=0$ 
established in Theorem~\juergen{4.1} is uniquely
determined. 
\Ethm
\Bdim
\pier{We point out that the third condition in \eqref{embedABC} is a straightforward consequence of the first and second ones. The continuity of the embedding $\VBs \subset \VAr$} \gianni{implies the existence of a constant~$\kappa$
(which we will refer~to) such that
\Beq
  \|\Ar v\|^2\,\le\,\kappa\bigl(\|\Bs v\|^2+\|v\|^2\bigr)\quad\mbox{for all }\,\,v\in \VBs.
  \label{compa}
\Eeq
}%
Let $(\mu_i,\phi_i,S_i)$, $i=1,2$, be two solution triples. We denote $w_i:=\juergen{\alpha}\mu_i+\phi_i$, for $i=1,2$, and
set $\mu:=\mu_1-\mu_2$, $\phi:=\phi_1-\phi_2$, $S:=S_1-S_2$, and $w:=w_1-w_2$. Then we have a.e. in $(0,t)$, and for $i=1,2$, that
\begin{align}
\label{mini1}
&\langle \dt w_i,v\rangle_{\VAr}+(\Ar\mu_i,v)\,=\,(P(\phi_i)(S_i-\mu_i),v)\quad\forall\,v\in\VAr,\\[1mm]
\label{mini2}
&\alpha(\Bs\phi_i,\Bs(\phi_i-v))+\alpha\iO F_1(\phi_i)+\alpha(f_2(\phi_i),\phi_i-v)+(\phi_i,\phi_i-v)\non\\[0.5mm]
&\le\,(w_i,\phi_i-v)+\alpha\iO F_1(v)\quad\forall\,v\in\VBs,\\[1mm]
\label{mini3}
&\langle \dt S_i,v\rangle_{\VCt}+(\Ct S_i,\Ct v)\,=\,(P(\phi_i)(\mu_i-S_i),v)\quad\forall\,v\in\VCt. 
\end{align}
Next, we insert $v=\phi_2$ in the inequality \eqref{mini2} for $i=1$, $v=\juergen{\phi_1}$ in the inequality for $i=2$,
add the resulting inequalities and multiply the result by a positive constant $M$ which is yet to be   
specified. Then we integrate over $(0,t)$\pier{,} where $t\in (0,T)$. Note that all of the terms involving $F_1$ cancel. 
Hence, also using \eqref{hpFdue}, we obtain the 
inequality
\begin{align*}
&M\alpha\int_{Q_t}|\Bs\phi|^2\,+\,M(1-\alpha L)\int_{Q_t}|\phi|^2\,\le\,M\int_{Q_t}w\,\phi\,,
\end{align*}
and Young's inequality yields that for every $\delta>0$ (which is yet to be chosen) it holds that
\begin{align}
\label{mini4}
M\alpha\int_{Q_t}|\Bs\phi|^2\,+\,(M(1-\alpha L)-\delta)\int_{Q_t}|\phi|^2\,\le\,\frac{M^2}
{4\delta}\int_{Q_t}|w|^2\,.
\end{align}
Now we subtract the equations \eqref{mini1} for $i=1,2$ from each other and insert $v=w$ in the resulting equation. Similarly,
we subtract the equations \eqref{mini3} for $i=1,2$ from each other and insert $v=S$ in the resulting equation. 
Finally, we add the two results. Integration over $(0,t)$ then yields the identity
\begin{align}
\label{mini5}
&\frac 12\,\|w(t)\|^2\,+\,\alpha\int_{Q_t}|\Ar\mu|^2\,+\,\int_{Q_t}\Ar \mu\,\Ar\phi\,+\,\frac 12\,\|S(t)\|^2
\,+\int_{Q_t}|\Ct S|^2\non\\
&=\,\int_{Q_t} \bigl(P(\phi_1)(S_1-\mu_1)-P(\phi_2)(S_2-\mu_2)\bigr)\,(w-S)\,.
\end{align}
Now observe that Young's inequality and \eqref{compa} yield that
\begin{align}
\label{mini6}
&-\int_{Q_t}\Ar\mu \,\Ar\phi\,\le\,\frac \alpha 2\int_{Q_t}|\Ar\mu|^2\,+\,\frac 2 \alpha\int_{Q_t}|\Ar\phi|^2\non\\
&\le\,\frac \alpha 2\int_{Q_t}|\Ar\mu|^2\,+\,\frac {2\kappa} \alpha\int_{Q_t}(|\Bs\phi|^2+|\phi|^2)\,.
\end{align}
It remains to estimate the \rhs\ of \eqref{mini5} which we denote by $Z$. We have
\begin{align}
\label{mini7}
Z\,&=\int_{Q_t} \bigl(P(\phi_1)-P(\phi_2))(S_1-\mu_1)(w-S)\,+\int_{Q_t}P(\phi_2)(S-\mu)\,(w-S)\non\\[1mm]
&=: \,Z_1+Z_2,
\end{align}
with obvious notation. Using the H\"older and Young inequalities, and invoking \eqref{embedABC}, we 
see that
\begin{align}
\label{mini8}
|Z_1|\,&\le\,C_1\int_0^t\|\phi(s)\|_{L^4(\Omega)}\,\bigl(\|S_1(s)\|_{L^4(\Omega)}
+\|\mu_1(s)\|_{L^4(\Omega)}\bigr)\,\bigl(\|w(s)\|+\|S(s)\|\bigr)\,ds\non\\
&\le\,\delta\int_0^t\|\phi(s)\|^2_{\VBs}\,ds\,+\,\frac{C_2}\delta\int_0^t
\Phi(s)\,\bigl(\|w(s)\|^2+\|S(s)\|^2\bigr)\,ds\,,
\end{align}
where the function $\,\,\Phi(s):=\|S_1(s)\|^2_{\VCt}\,+\,\|\mu_1(s)\|_{\VAr}^{\gianni 2}\,$\, is known to belong
to $L^1(0,T)$. Here, and in the remainder of the proof, $C_i$, $i\in\enne$, denote positive
constants that depend only on the global data of the system. 

Finally, we estimate $Z_2$. Omitting an obvious nonpositive term, we have, by virtue of Young's inequality,
and since $\,\alpha\mu=w-\phi$,
\begin{align}
\label{mini9}
|Z_2|\,&\le\,C_3\int_0^t\bigl(\|S(s)\|\,\|w(s)\|\,+\,\|S(s)\|\,\|\mu(s)\|\,+\,\|\mu(s)\|\,\|w(s)\|\bigr)
\,ds\non\\
&\le \,\delta\,\alpha^2\int_{Q_t}|\mu|^2\,+\,C_4\Big(1+ \,\frac 1{\delta\alpha^2}\Big)
\int_{Q_t}\bigl(|S|^2+|w|^2\bigr)
\non\\
&\le\,2\delta\int_0^t\|\phi(s)\|^2_{\VBs}\,ds \,+\, C_5\Big(1+\delta+\,\frac 1{\delta\alpha^2}\big)
\int_{Q_t}\bigl(|S|^2+|w|^2\Big)\,.
\end{align}
Combining \eqref{mini4}--\eqref{mini9}, we have thus shown the estimate
\begin{align}
\label{mini10}
&\Big(M\alpha-\,\frac{2\kappa}\alpha\,-3\delta\Big)\int_{Q_t}|\Bs\phi|^2\,+\,\Big(M(1-\alpha L)-\,
\frac{2\kappa}\alpha\,-4\delta\Big)\int_{Q_t}|\phi|^2\non\\
&+\,\frac 12\bigl(\|w(t)\|^2+\|S(t)\|^2\bigr)\,+\,\frac\alpha 2\int_{Q_t}|\Ar\mu|^2
\,+\,\int_{Q_t}|\Ct S|^2\non\\
&\le\int_0^t\Big[\frac{C_2}\delta\,\Phi(s)\,+\,C_6\,\juergen{\bigl(M^2+1\bigr)}\Big(1+\delta+\,\frac 1\delta\,+\,
\frac 1{\delta\alpha^2}\Big)\Big]
\,\bigl(\|w(s)\|^2+\|S(s)\|^2\bigr)\,ds\,.
\end{align}
At this point, we make the choices 
$$
M>M_0:=\max\,\left\{\frac {2\kappa}{\alpha^2}\,,\,\frac{2\kappa}{\alpha(1-\alpha L)}\right\}\quad
\mbox{and}\quad 0<\delta<\frac 14(M-M_0)\,\min\,\{\alpha,1-\alpha L\}\,.
$$
Then the brackets in the first two terms on the \lhs\ become positive, and we may apply Gronwall's
lemma to conclude that $w=S=\phi=0$, whence also $\mu=0$.
\Edim

\vspace{3mm}
\Brem
It ought to be clear from the above arguments that in the case that controls $u_\mu$, $u_\phi$, $u_S$ 
in $L^2(0,T;H)$ are added to the \rhs s of \eqref{prima}--\eqref{terza}, we have an existence 
result resembling Theorem~\gianni{\ref{Atozero}}, 
and, under the assumptions of Theorem~\gianni{\ref{UniqAzero}}, we obtain a corresponding continuous dependence
result in the norms appearing on the \lhs\ of \eqref{mini10}. 
\Erem

\section{The case $\mathbf{\alpha\searrow0,\,\beta\searrow0.}$}
\setcounter{equation}{0}

In this section, we investigate the asymptotic behavior of the solutions $(\muab,\phiab,\Sab)$
as $\alpha\searrow0$ and $\beta\searrow0$. Quite unexpectedly, in this case the additional assumption
{\bf (A6)} is not needed. In a sense, this means that the presence of a strong perturbation
$\,\alpha\dt\mu$ as in the previous section does not just produce an approximation but really changes
the character of the unperturbed system if $\alpha$ is too large. On the other hand, we have to assume:

\vspace{2mm}\noindent
{\bf (A7)} \,\,The eigenvalue $\lambda_1$ is positive.

\vspace{2mm}\noindent
Recall that then the mapping $\,v\mapsto \|\Ar v\|\,$ defines a norm on $\VAr$ which is equivalent to
the graph norm $\norma\cpto_{\VAr}$ \gianni{(see Remark~\ref{Meanvalue})}.

To begin with, let $\{\alpha_n\}_{n\in\enne}\subset (0,1]$ and $\{\beta_n\}_{n\in\enne}\subset (0,1]$
be sequences such that $\alpha_n\searrow0$ and $\beta_n\searrow0$, and let $(\mu_{\alpha_n,\beta_n},
\phi_{\alpha_n,\beta_n},S_{\alpha_n,\beta_n})$ denote solutions to \eqref{prima}--\eqref{cauchy} in the sense of
Theorem~\gianni{\ref{Wellposed}} associated with $(\alpha,\beta)=(\alpha_n,\beta_n)$, for $n\in\enne$. 
According to \eqref{bounds1},
and invoking {\bf (A7)}, we may without loss of generality assume that there are limits $\,\zeta,
\mu_{0,0},\phi_{0,0},S_{0,0}\,$ such that,
as $n\to\infty$,
\begin{align}
\label{convc1}
&\alpha_n\,\mu_{\alpha_n,\beta_n}\,\to\,0\quad\mbox{strongly in }\,L^\infty(0,T;H),\\[1mm]
\label{convc2}
&\beta_n\,\dt\phi_{\alpha_n,\beta_n}\,\to\,0\quad\mbox{strongly in }\,L^2(0,T;H),\\[1mm]
\label{convc3}
&\dt(\alpha_n\,\mu_{\alpha_n,\beta_n}+\phi_{\alpha_n,\beta_n})\,\to\,\zeta \quad\mbox{weakly in }\,L^2(0,T;V_A^{-\rho}),
\\[1mm]
\label{convc4}
&\mu_{\alpha_n,\beta_n}\,\to\,\mu_{0,0}\quad\mbox{weakly in }\,L^2(0,T;\VAr),\\[1mm]
\label{convc5}
&\phi_{\alpha_n,\beta_n}\,\to\,\phi_{0,0}\quad\mbox{weakly-star in }\,L^\infty(0,T;\VBs),\\[1mm]
\label{convc6}
&S_{\alpha_n,\beta_n}\,\to\,S_{0,0} \quad\mbox{weakly-star in }\,H^1(0,T;V_C^{-\tau})\cap L^\infty(0,T;H)\cap
L^2(0,T;\VCt).  
\end{align}
From \pier{\eqref{convc1}, \eqref{convc3} and \eqref{convc5} it follows that 
$\,\zeta=\dt\phi_{0,0}$ and, in addition,
\Beq
  \alpha_n\,\mu_{\alpha_n,\beta_n}+\phi_{\alpha_n,\beta_n} \to \phi_{0,0}
  \quad \hbox{weakly in $\H1{\VA{-\rho}}$},
  \label{convc7}
\Eeq
whence also
\Beq
  (\alpha_n\,\mu_{\alpha_n,\beta_n}+\phi_{\alpha_n,\beta_n})(0) \to \phi_{0,0}(0)
  \quad \hbox{weakly in $\VA{-\rho}$}.
  \label{convc8}
\Eeq
Then, in view of \eqref{convc6}--\eqref{convc8}, it turns out that both the initial conditions in} \eqref{cauchyazero} are fulfilled.

Next, we observe that we can argue exactly 
\gianni{as we did in the previous section to obtain~\eqref{conva6}.
Hence, we}
infer that the sequence $\,\{\alpha_n\,\mu_{\alpha_n,\beta_n}+\phi_{\alpha_n,\beta_n}\}_{n\in\enne}\,$
converges strongly in $L^2(0,T;H)$. 
We thus find from \eqref{convc1} that\pier{%
\begin{align}
\label{convc10}
&\phi_{\alpha_n,\beta_n}\,\to\,\phi_{0,0}\quad\mbox{strongly in \,$L^2(0,T;H)\,$ and}
\non\\
& \quad \pier{\phi_{\alpha_n,\beta_n}(t) \,\to\,\phi_{0,0}(t) \quad\mbox{strongly in \,$H$,
 for a.e.\ $t\in (0,T)$,}}
\end{align}
}%
the latter without loss of generality. 
Consequently, by Lipschitz continuity, \pier{we have that}
\begin{equation}
\label{convc11}
f_2(\phi_{\alpha_n,\beta_n})\,\to\,f_2(\phi_{0,0})\quad \pier{\mbox{and}} 
\quad P(\phi_{\alpha_n,\beta_n})\,\to\,P(\phi_{0,0})\quad\mbox{strongly in }\,L^2(0,T;H).
\end{equation}
\gianni{From this point, we may follow the lines of the previous sections to conclude the following result}.

\Bthm
\label{ABtozero}
Assume that {\bf (A1)}--{\bf (A4)} and {\bf (A7)} are fulfilled
\gianni{and let the sequences $\{\alpha_n\}_{n\in\enne}\subset(0,1]$ and $\{\beta_n\}_{n\in\enne}\subset (0,1]$ 
satisfy $\alpha_n\searrow0$ and $\beta_n\searrow0$}.
Moreover, let $(\mu_{\alpha_n,\beta_n},\phi_{\alpha_n,\beta_n},S_{\alpha_n,\beta_n})$
be solutions \gianni{to the system \eqref{prima}--\eqref{cauchy}} in the sense of Theorem~\gianni{\ref{Wellposed}} 
for $(\alpha,\beta)=(\alpha_n,\beta_n)$ for $n\in\enne$. 
Then, there are a subsequence $\{n_k\}_{k\in\enne}$ of $\enne$ and 
a triple $(\mu_{0,0}, \phi_{0,0},S_{0,0})$ such that the following holds true:
\begin{align}
\label{cokc1}
&\alpha_{n_k}\,\mu_{\alpha_{n_k},\beta_{n_k}}\,\to\,0\quad\mbox{strongly in }\,L^\infty(0,T;H),\\[1mm]
\label{cokc2}
&\beta_{n_k}\,\dt\phi_{\alpha_{n_k},\beta_{n_k}}\,\to\,0\quad\mbox{strongly in }\,L^2(0,T;H),\\[1mm]
\label{cokc3}
&\dt(\alpha_{n_k}\,\mu_{\alpha_{n_k},\beta_{n_k}}+
\phi_{\alpha_{n_k},\beta_{n_k}})\,\to\,\dt\phi_{0,0} \quad\mbox{weakly in }\,L^2(0,T;V_A^{-\rho}),
\\[1mm]
\label{cokc4}
&\mu_{\alpha_{n_k},\beta_{n_k}}\,\to\,\mu_{0,0}\quad\mbox{weakly in }\,L^2(0,T;\VAr),\\[1mm]
\label{cokc5}
&\phi_{\alpha_{n_k},\beta_{n_k}}\,\to\,\phi_{0,0}\quad\mbox{weakly\juergen{-star} in }\,L^\infty(0,T;\VBs)
\,\mbox{ and strongly in }\,\,L^2(0,T;H),\\[1mm]
\label{cokc6}
&S_{\alpha_{n_k},\beta_{n_k}}\,\to\,S_{0,0} \quad\mbox{weakly-star in }\,H^1(0,T;V_C^{-\tau})\cap L^\infty(0,T;H)\cap
L^2(0,T;\VCt).  
\end{align}
Moreover, \gianni{$F_1(\phi_{0,0})\in\pier{\L\infty\Luno}$, and} $(\mu_{0,0}, \phi_{0,0},S_{0,0})$\, 
is a solution to \eqref{prima}--\eqref{terza} for $\alpha=\beta=0$ \gianni{that}
satisfies the initial \gianni{conditions \eqref{cauchyazero}}. 
\Ethm 

It seems difficult to prove a uniqueness result for the solution to the limiting problem
under rather general assumptions.
The arguments developed in \cite{CGH} and \cite{FGR} strongly use the fact that 
$A^{2\rho}$, $B^{2\sigma}$ and $C^{2\tau}$ are the same operator 
(namely, the Laplace operator with zero Neumann boundary conditions),
and this kind of assumption is quite unpleasant in the \juergen{context} of the present paper.
Thus, we prefer to keep the operators $A$, $B$ and $C$ and the exponents $\rho$, $\sigma$ and $\tau$
independent from each other.
To do this, we have to make restrictive assumptions and compatibility conditions on the structure of the system.
We recall that $\lambda'_j$ and $L$ are the eigenvalues of $B$ and the \Lip\ constant of~$f_2$, respectively,
and \juergen{make the following requirement:}

\medskip\noindent
{\bf (A8)}\quad
The following conditions are satisfied:
\begin{align}
  & F_1 \in C^1(\erre), \quad
  |f_1(s)| \leq \hat c_5 (|s|^3 + 1)
  \aand
  (f_1(s)-f_1(s'))(s-s') \geq \gamma \, |s-s'|^2 
  \qquad
  \non
  \\
  & \hbox{for some constants $\hat c_5\,,\,\gamma>0$ and every $s,s'\in\erre$}\,.
  \label{hpf1}
  \\
  &\mbox{We have }\, (\lambda'_1)^{2\sigma} + \gamma > L \,.
  \label{compat}
  \\
  & \hbox{It holds the continuous embedding} \quad
  \VB\sigma \subset \Lx 4 \,.
  \label{embB}
  \\
  & \hbox{$P$ is a positive constant $P_0$}\,.
  \label{hpPconst}
\end{align}

\Brem
\label{RemA8}
The first condition excludes singular potentials and holds for the regular potential $F_{reg}$ given by~\eqref{regpot}.
As for the strong monotonicity condition in~\eqref{hpf1},
one can split $F$ into $F_1+F_2$ according to \accorpa{hpFuno}{hpbelow} by setting
$$
  F_1(s) := \frac 14 \, s^4 + \frac 12 \, s^2
  \aand
  F_2(s) := - s^2 + \frac 14 \,.
$$
Then, $f_1'(s)=F_1''(s)=3s^2+1\geq1$ for every $s\in\erre$,
so that we can take $\gamma=1$ in~\eqref{hpf1}.
The compatibility condition \eqref{compat} is not satisfied by double-well potentials if $\lambda'_1=0$,
since it reduces to $\gamma>L$ in this case and is satisfied \juergen{only} if $F$ is convex.
On the contrary, if we consider the above splitting of the regular potential~$F_{reg}$, \juergen{then}
we see that \eqref{compat} holds if $\lambda'_1$ is large enough, 
namely if $\lambda'_1>1$, since $\gamma=1$ and $L=2$.
If $B$ is, e.g., the Laplace operator with zero Dirichlet boundary conditions,
\juergen{then} this kind of assumption on $\lambda'_1$ is satisfied provided that $\Omega$ is small enough
\juergen{within} a class of domains having the same shape.
Finally, embeddings similar to \eqref{embB} have already been commented in Remark~\ref{RemembedAC},
and \eqref{hpPconst} is not realistic in the framework of the tumor model, unfortunately.
\Erem

\Bthm
\label{UniqABzero}
Besides {\bf(A1)}--{\bf(A4)}, assume that {\bf(A8)} is satisfied as well.
Then the problem \accorpa{prima}{cauchy} with $\alpha=\beta=0$
has at most one solution satisfying \pier{\accorpa{regmu}{L1Q}}.
\Ethm

\Bdim
First of all, we prove that any solution $(\mu,\phi,S)$ satisfies an equation like~\eqref{Iseconda}.
The inequality \eqref{seconda} with $\beta=0$ becomes
\Beq
  \bigl( B^\sigma \phi(t) , B^\sigma( \phi(t)-w) \bigr)
  + \iO F_1(\phi(t))
  \leq \bigl( \mu(t) - f_2(\phi(t)) , \phi(t)-w \bigr)
  + \iO F_1(w)
  \non
\Eeq
\aat\ and every $w\in\VB\sigma$.
We can take, in particular, $w=\phi(t)-\eps v$ with any $v\in\VB\sigma$ and $\eps\in(0,1)$.
We then obtain that
\Beq
  \bigl( B^\sigma \phi(t) , B^\sigma v \bigr)
  + \iO \frac {F_1(\phi(t)) - F_1(\phi(t)-\eps v)} \eps
  \leq \bigl( \mu(t) - f_2(\phi(t)) , v \bigr).
  \non
\Eeq
Now, for fixed~$t$, we apply the mean value theorem, the growth condition in~\eqref{hpf1} and the Young inequality.  \pier{Using the integral remainder,}
we have \aeO\ that 
\Bsist
  && |F_1(\phi(t)) - F_1(\phi(t)-\eps v)|
  = \pier{\biggl| \int_0^1 f_1 (\phi(t)-\badtheta\eps v)\eps v \, d \badtheta \biggr|} 
  \non
  \\
  && \leq \pier{ \hat c_5 \bigl( (|\phi(t)| + \eps |v|)^3 + 1 \bigr)  \eps \, |v|} 
  \leq \eps \, C (|\phi(t)|^4 + |v|^4\juergen{{}+1} )\,,
  \non
\Esist
\pier{with} some constant $C$ proportional to~$\hat c_5$.
Since $|\phi(t)|^4+|w|^4\in\Luno$ thanks to the embedding~\eqref{embB},
we can apply the Lebesgue dominated convergence theorem 
and let $\eps$ tend to zero. 
We deduce an inequality holding for every $v\in\VB\sigma$.
Since $v$ is arbitrary, we obtain the equality
\Beq
  \bigl( B^\sigma \phi(t) , B^\sigma v \bigr)
  + \iO f_1(\phi(t)) \, v
  = \bigl( \mu(t) - f_2(\phi(t)) , v \bigr),
  \label{secondaeq}
\Eeq
which \juergen{is valid} \aat\ and every $v\in\VB\sigma$.

We notice that the integral in \eqref{secondaeq} cannot be written as $(f_1(\phi(t)),v)$,
since we do not know that $f_1(\phi(t))$ belongs to~$H$.
However, it belongs to $\Lx{4/3}$, 
since $\phi(t)\in\VB\sigma\subset\Lx4$ and the growth condition in~\eqref{hpf1}) is in force.
For this reason, if we also assume that $\VB\sigma$ is dense in~$\Lx4$,
\juergen{then} we can write \eqref{Iseconda} in the sense of $\VB{-\sigma}$
by accounting for the consequent embedding $\Lx{4/3}\subset\VB{-\sigma}$.
However, we will use just \eqref{secondaeq} as it~is.

We are ready to prove uniqueness.
We pick two solutions $(\mu_i,\phi_i,S_i)$, $i=1,2$, and set for convenience
$\mu:=\mu_1-\mu_2$, $\phi_1-\phi_2$, and $S:=S_1-S_2$.
We write \eqref{prima} and \eqref{terza} for both solutions and take the differences.
By recalling \eqref{hpPconst}, we \pier{arrive at} the identities
\begin{align}
  & \pier{\< \,\dt \phi(t) , v \,>_\VAr +{}} \bigl( \Ar \mu(t) , \Ar v \bigr)
  \,=\, P_0 \bigl( S(t) - \mu(t) , v \bigr),
  \non
  \\
  & \< \,\dt S(t) , v \,>_\VCt
  \,+\, \bigl( \Ct S(t) ,C^\tau v \bigr)
  \,=\, - P_0 \bigl( S(t) - \mu(t) , v \bigr),
  \non
\end{align}
which hold for every $v\in\VA\rho$ and $v\in\VC\tau$, respectively.
Now, we integrate these equations with respect to time \pier{and, 
for $X$ Banach space and $w\in\L1X$, we use the notation
\Beq
  (1*w)(t) := \iot w(s) \ds
  \quad \hbox{for every $t\in[0,T]$}. 
  \non
\Eeq
Hence, we obtain}
\begin{align}
  & \pier{\bigl( \phi(t) , v \bigr) +{}} \bigl( \Ar (1*\mu)(t) , \Ar v \bigr)
  \,=\, P_0 \bigl( (1*(S-\mu))(t) , v \bigr),
  \non
  \\
  & \bigl( S(t) , v \bigr)
  \,+\, \bigl( \Ct (1*S)(t) ,C^\tau v \bigr)
  \,=\, - P_0 \bigl( (1*(S-\mu))(t) , v \bigr),
  \non
\end{align}
for the same test functions $v$ as before.
At this point, we choose $v=\mu(t)$ and $v=S(t)$ in these identities, respectively.
Then, we integrate with respect to time, sum up and rearrange.
We \pier{deduce}, for every $t\in[0,T]$,~that
\begin{align}
  & \frac 12 \, \norma{\Ar (1*\mu)(t)}^2
  + \iot \norma{S(s)}^2 \ds
  \non \\
  &+ \frac 12 \, \norma{\Ct (1*S)(t)}^2
  + \frac {P_0} 2 \, \norma{1*(S-\mu))(t)}^2
  = \pier{{}- \iot \bigl( \phi(s) , \mu(s) \bigr)\ds{}} \pcol{.}
  \label{pier1}
\end{align}
\pier{Next,} we write equation \eqref{secondaeq} for both solutions
and choose $v=\phi(t)$ in the difference. \pier{By integrating with respect to time, we infer that 
\begin{align}
  &\iot \norma{B^\sigma \phi(s)}^2\ds
  + \int_{Q_t} \bigl( f_1(\phi_1) - f_1(\phi_2) \bigr) \phi
  \non \\
  &= \iot \bigl( \mu(s) , \phi(s) \bigr)\ds  
  - \int_{Q_t} \bigl( f_2(\phi_1) - f_2(\phi_2) \bigr) \phi.
  \non
\end{align}
Now, we account for the obvious inequality $\norma{B^\sigma v}^2\geq(\lambda'_1)^{2\sigma}\norma v^2$,
the assumption~\eqref{hpf1}, and the \Lip\ continuity of~$f_2$ \pier{(cf.~\eqref{hpFdue} and \eqref{deffunodue}), in order to deduce that}
\Beq
  (\lambda'_1)^{2\sigma} \iot \norma{\phi(s)}^2\ds 
  + \gamma \int_{Q_t} |\phi|^2
  \leq \iot \bigl( \mu(s) , \phi(s) \bigr)\ds  
  + L  \int_{Q_t} | \phi |^2
  \label{pier2}
\Eeq
for every $t\in[0,T]$. Now, we add \eqref{pier1} and \eqref{pier2}, note that there is 
a cancellation and finally apply \eqref{compat}. Hence, we conclude in particular that $\phi=0$, $S=0$, $1*(S-\mu)=0$. The latter implies that $S-\mu=0$, whence $\mu=0$ as well. 
We have thus proved that $\phi_1=\phi_2$, $S_1=S_2$ and $\mu_1=\mu_2$.}
\Edim

\Brem
\label{RemAposBzero}
The same assumption {\bf(A8)} (possibly reinforced by also supposing that $f_1$ and $f_2$ are $C^1$ functions)
can be used to prove a regularity result in the case $\alpha>0$ and $\beta=0$.
Here we sketch a formal proof under suitable assumptions on the initial data,
by observing that {\bf(A8)} ensures the validity of \eqref{secondaeq} also in this case.
Indeed, the argument used in the proof of Theorem~\ref{UniqABzero} to derive \eqref{secondaeq}
only regards the variational inequality satisfied by $\phi$ for $\beta=0$ and thus still holds if $\alpha>0$.
We differentiate \eqref{prima} and \eqref{terza} with respect to time
and test the equalities we get by $\dt\mu$ and~$\dt S$, respectively.
At the same time, we test \pier{the time derivative of} \eqref{secondaeq} by~$\dt\phi$.
Then, we sum up and integrate over~$(0,t)$.
The terms involving the product $(\dt\mu,\dt\phi)$ cancel each other, and we obtain 
(by~omitting the integration variable $s$ to shorten the lines) the identity
\Bsist
  && \alpha \iot \norma{\dt\mu}^2 \ds
  + \frac 12 \, \norma{A^\rho\mu(t)}^2
  + \iot \norma{B^\sigma\dt\phi}^2 \ds
  + \intQt f_1'(\phi) \, |\dt\phi|^2 
  \non
  \\
  && \quad {}
  + \iot \norma{\dt S}^2 \ds
  + \frac 12 \, \norma{C^\tau S(t)}^2 
  + \frac {P_0} 2 \, \norma{S(t)-\mu(t)}^2
  \non
  \\
  && = \frac 12 \, \norma{A^\rho\mu_0}^2
  + \frac 12 \, \norma{C^\tau S_0}^2 
  + \frac {P_0} 2 \, \norma{S_0-\mu_0}^2
  - \intQt f_2'(\phi) \, |\dt\phi|^2 \,.
  \non
\Esist
Now, from one side, we have that $\norma{B^\sigma\dt\phi}^2\geq\juergen{(\lambda_1')}^{2\sigma}\norma{\dt\phi}^2$.
On the other hand, \eqref{hpf1} and \eqref{hpFdue} imply that 
$\,f_1'(\phi)\geq\gamma\,$ and $\,|f_2'(\phi)|\leq L\,$ \aeQ.
Therefore, we derive, for every $\delta>0$,~that
\Bsist
  && \iot \norma{B^\sigma\dt\phi}^2 \ds
  + \intQt f_1'(\phi) \, |\dt\phi|^2 
  + \intQt f_2'(\phi) \, |\dt\phi|^2 
  \non
  \\
  && \pier{{}\ge{}} \delta \iot \norma{B^\sigma\dt\phi}^2 \ds
  + \bigl( (1-\delta) \juergen{(\lambda_1')}^{2\sigma} + \gamma - L \bigr) \iot \norma{\dt\phi}^2 \ds \,.
  \non
\Esist
By choosing $\delta>0$ such  that $(1-\delta)\juergen{(\lambda_1')}^{2\sigma}+\gamma>L$ on account of~\eqref{compat}, 
we conclude that
\Bsist
  && \norma{\dt\mu}_{\L2H}
  + \norma{A^\rho\mu}_{\L\infty{\VA\rho}}
  + \norma{\dt\phi}_{\L2{\VB\sigma}}
  \non
  \\
  && \quad {}
  + \norma{\dt S}_{\L2H}
  + \norma{C^\tau S}_{\L\infty H}
  + \norma{S-\mu}_{\L\infty H}
  \leq C\,,
  \non
\Esist
where $\,C\,$ depends \juerg{only} 
on the structural assumptions and the norms of the initial data involved in the calculation.
\Erem



\section*{Acknowledgments}
\pier{This research was supported by the Italian Ministry of Education, 
University and Research~(MIUR): Dipartimenti di Eccellenza Program (2018--2022) 
-- Dept.~of Mathematics ``F.~Casorati'', University of Pavia. 
In addition, PC and CG gratefully acknowledge some other 
financial support from the GNAMPA (Gruppo Nazionale per l'Analisi Matematica, 
la Probabilit\`a e le loro Applicazioni) of INdAM (Isti\-tuto 
Nazionale di Alta Matematica).}

\footnotesize


{


\vspace{3truemm}

\Begin{thebibliography}{10}


\pcol{\bibitem{BKM}
B. Baeumer, M. Kov\'acs, M. M. Meerschaert, 
Numerical solutions for fractional reaction-diffusion equations,
{\it Comput. Math. Appl.}  {\bf 55}  (2008), 2212--2226.}

\bibitem{BLM}
N. Bellomo, N. K. Li, P. K. Maini,
On the foundations of cancer modelling: Selected topics, speculations, and perspectives,
\textit{Math. Models Methods Appl. Sci.} \textbf{18} (2008), 593--646.

\bibitem{BCG}
S. Bosia, M. Conti, M. Grasselli,
On the Cahn--Hilliard--Brinkman system,
\textit{Commun. Math. Sci.} \textbf{13} (2015), 1541--1567.

\bibitem{Brezis}
H. Brezis,
``Op\'erateurs maximaux monotones et semi-groupes de contractions
dans les espaces de Hilbert'',
North-Holland Math. Stud.
{\bf 5},
North-Holland,
Amsterdam,
1973.

\bibitem{Cahn}
J. W. Cahn, J. E. Hilliard,
Free energy of a nonuniform system. I. Interfacial free energy,
{\it J. Chem. Phys.} {\bf 28} (1958), 258--267.

\bibitem{CRW}
C. Cavaterra, E. Rocca, H. Wu,
Long-time dynamics and optimal control of a diffuse interface model for tumor growth,
\textit{Appl. Math. Optim.} (2019), Online First 15 March 2019, 
https://doi.org/10.1007/s00245-019-09562-5.

\pcol{\bibitem{CB}
S. K. Chandra, M. K. Bajpai, 
Mesh free alternate directional implicit method based three dimensional 
super-diffusive model for benign brain tumor segmentation,
{\it Comput. Math. Appl.}  {\bf 77}  (2019), 3212--3223.}

\bibitem{CWSL}
Y. Chen, S. M. Wise, V. B. Shenoy, J. S. Lowengrub,
A stable scheme for a nonlinear multiphase tumor growth model with an elastic membrane,
\textit{Int. J. Numer. Methods Biomed. Eng.} \textbf{30} (2014), 726--754.

\bibitem{CG1} P. Colli, G. Gilardi, 
Well-posedness, regularity and asymptotic analyses 
for a fractional phase field system,
\pier{{\em Asymptot. Anal.} \textbf{114} (2019), 93--128.}

\bibitem{CGH}
P. Colli, G. Gilardi, D. Hilhorst,
On a Cahn--Hilliard type phase field system related to tumor growth,
\textit{Discrete Contin. Dyn. Syst.} \textbf{35} (2015), 2423--2442.

\bibitem{CGMR}
P. Colli, G. Gilardi, G. Marinoschi, E. Rocca, 
Sliding mode control for a phase field system related to tumor growth,
{\em Appl. Math. Optim.} \textbf{79} (2019), 647--670.

\bibitem{CGRS}
P. Colli, G. Gilardi, E. Rocca, J. Sprekels,
Vanishing viscosities and error estimate for a Cahn--Hilliard type phase field system related to tumor growth,
\textit{Nonlinear Anal. Real World Appl.} \textbf{26} (2015), 93--108.

\pier{\bibitem{CGRS2}
P. Colli, G. Gilardi, E. Rocca, J. Sprekels, 
Asymptotic analyses and error estimates for a Cahn--Hilliard type phase field system modelling tumor growth,
\textit{Discrete Contin. Dyn. Syst. Ser. S.} {\bf  10} (2017), 37--54. 
\bibitem{CGRS1}
P. Colli, G. Gilardi, E. Rocca, J. Sprekels,
Optimal distributed control of a diffuse interface model of tumor growth,
\textit{Nonlinearity} \textbf{30} (2017), 2518--2546.}

\bibitem{CGS18}
P. Colli, G. Gilardi, J. Sprekels,
Well-posedness and regularity for a generalized fractional Cahn--Hilliard system,
{\em Atti Accad. Naz. Lincei Rend. Lincei Mat. Appl.}
\pier{{\bf 30} (2019) 437--478.}

\bibitem{CGS19}
P. Colli, G. Gilardi, J. Sprekels,
Optimal distributed control of a generalized fractional Cahn--Hilliard system,
{\em Appl. Math. Optim.} (2018), Online First 15 November 2018, 
https://doi.org/10.1007/s00245-018-9540-7.

\bibitem{CGS21}
P. Colli, G. Gilardi, J. Sprekels,
Deep quench approximation and optimal control of 
general Cahn--Hilliard systems with fractional
operators and double-obstacle potentials, 
\pier{\textit{Discrete Contin. Dyn. Syst. Ser. S.}, to appear
(see also preprint arXiv:1812.01675 [math.AP] (2018), pp. 1--32).}


\bibitem{CGS22}
P. Colli, G. Gilardi, J. Sprekels,
Longtime behavior for a generalized Cahn--Hilliard system with fractional operators,
preprint arXiv:1904.00931 [math.AP] (2019), pp. 1--18, submitted.

\bibitem{CGS23}
P. Colli, G. Gilardi, J. Sprekels,
Well-posedness and regularity for a fractional tumor growth model, 
{\em Adv. Math. Sci. Appl.} \pier{{\bf 28} (2019) 343--375.}

\bibitem{CGS25}
P. Colli, G. Gilardi, J. Sprekels,
A distributed control problem for a fractional tumor growth model,
\pcol{\textit{Mathematics}, to appear
(see also preprint arXiv:1907.10452 [math.OC] (2019), pp. 1--35).}

\bibitem{ConGio}
M. Conti, A. Giorgini,
The three-dimensional Cahn--Hilliard--Brinkman system with unmatched viscosities,
preprint hal-01559179 (2018), pp.~1--34.

\bibitem{CLLW}
V. Cristini, X. Li, J. S. Lowengrub, S. M. Wise, Nonlinear simulations of solid tumor growth using a 
mixture model: invasion and branching, \textit{J. Math. Biol.} {\bf 58} (2009), 723--763.  

\bibitem{CL2010}
V. Cristini, J. S. Lowengrub,
``Multiscale Modeling of Cancer: An Integrated Experimental 
and Mathematical Modeling Approach'',
Cambridge Univ. Press, Cambridge, 2010.

\bibitem{DFRSS}
M. Dai, E. Feireisl, E. Rocca, G. Schimperna, M. Schonbek,
Analysis of a diffuse interface model for multi-species tumor growth,
\textit{Nonlinearity} \textbf{30} (2017), 1639--1658.

\bibitem{DG}
F. Della Porta, M. Grasselli,
On the nonlocal Cahn--Hilliard--Brinkman and Cahn--Hilliard--Hele--Shaw systems,
\textit{Commun. Pure Appl. Anal.} \textbf{15} (2016), 299--317,
Erratum: \textit{Commun. Pure Appl. Anal.} \textbf{16} (2017), 369--372.

\bibitem{DGG}
F. Della Porta, A. Giorgini, M. Grasselli,
The nonlocal Cahn--Hilliard--Hele--Shaw system with logarithmic potential,
\textit{Nonlinearity} \text{31} (2018), 4851--4881.

\bibitem{EGAR}
M. Ebenbeck, H. Garcke,
Analysis of a Cahn--Hilliard--Brinkman model for tumour growth with chemotaxis,
{\it J. Differential Equations} \textbf{266} (2019), 5998--6036.

\pcol{\bibitem{EGPS}
G. Estrada-Rodriguez, H. Gimperlein, K. J. Painter, J. Stocek,  
Space-time fractional diffusion in cell movement models with delay.
{\it Math. Models Methods Appl. Sci.}  {\bf 29} (2019), 65--88.}

\pcol{\bibitem{EL}
L. R. Evangelista, E. K. Lenzi, 
``Fractional diffusion equations and anomalous diffusion'',
Cambridge University Press, Cambridge, 2018.}

\bibitem{FBG2006}
A. Fasano, A. Bertuzzi, A. Gandolfi,
Mathematical modeling of tumour growth and treatment,
Complex Systems in Biomedicine, Springer, Milan, 2006, \pier{pp.~71--108.}

\bibitem{FLR}
S. Frigeri, K.~F. Lam, E. Rocca, 
On a diffuse interface model for tumour growth with non-local 
interactions and degenerate mobilities, 
in ``Solvability, regularity, and optimal control of boundary 
value problems for PDEs'', 
P.~Colli, A.~Favini, E.~Rocca, G.~Schimperna, J.~Sprekels~(ed.), 
Springer INdAM Series~{\bf 22}, Springer, Cham, 2017, pp.~217--254.

\bibitem{FW2012}
X. Feng, S. M. Wise,
Analysis of a Darcy--Cahn--Hilliard diffuse interface model for the Hele--Shaw flow and its fully discrete finite element approximation,
\textit{SIAM J. Numer. Anal.} \textbf{50} (2012), 1320--1343.

\bibitem{Fri2007}
A. Friedman,
Mathematical analysis and challenges arising from models of tumor growth,
\textit{Math. Models Methods Appl. Sci.} \textbf{17} (2007), 1751--1772.

\bibitem{Lowen10}
H. B. Frieboes, F. Jin, Y. L. Chuang, S. M. Wise, J. S. Lowengrub, V. Cristini,
Three-dimensional multispecies nonlinear tumor growth - II: tumor invasion and angiogenesis,
{\it J. Theoret. Biol.} {\bf 264} (2010), 1254--1278.

\bibitem{FGR}
S. Frigeri, M. Grasselli, E. Rocca,
On a diffuse interface model of tumor growth,
\textit{European J. Appl. Math.} \textbf{26} (2015), 215--243.

\bibitem{FLRS}
S. Frigeri, K. F. Lam, E. Rocca, G. Schimperna,
On a multi-species Cahn--Hilliard--Darcy tumor growth model with singular potentials,
\textit{Comm Math Sci.} \textbf{16} (2018), 821--856.

\bibitem{GL2016}
H. Garcke, K. F. Lam,
Global weak solutions and asymptotic limits of a Cahn--Hilliard--Darcy system modelling tumour growth,
\textit{AIMS Mathematics} \textbf{1} (2016), 318--360.

\bibitem{GL2017-1}
H. Garcke,  K.~F. Lam,
Analysis of a Cahn--Hilliard system with non--zero Dirichlet 
conditions modeling tumor growth with chemotaxis,
{\it Discrete Contin. Dyn. Syst.} {\bf 37} (2017), 4277--4308.
\bibitem{GL2017-2}
H. Garcke, K.~F. Lam,
Well-posedness of a Cahn--Hilliard system modelling tumour
growth with chemotaxis and active transport,
{\it European J. Appl. Math.} {\bf 28} (2017), 284--316.

\bibitem{GL2018}
H. Garcke,  K. F. Lam,
On a Cahn--Hilliard--Darcy system for tumour growth with solution dependent source terms,
In: Rocca E., Stefanelli U., Truskinovsky L., Visintin A. (eds),
``Trends in Applications of Mathematics to Mechanics'', Springer INdAM Series, Vol. \textbf{27}, Springer, 2018.

\bibitem{GLNS}
H. Garcke, K. F. Lam, R. N\"urnberg, E. Sitka,
A multiphase Cahn--Hilliard--Darcy model for tumour growth with necrosis,
\textit{Math. Models Methods Appl. Sci.} \textbf{28} (2018), 525--577.

\bibitem{GLR}
H. Garcke, K. F. Lam, E. Rocca,
Optimal control of treatment time in a diffuse interface model for tumour growth,
\textit{Appl. Math. Optim.} {\bf 78} \pier{(2018)}, 495--544.

\bibitem{GLSS}
H. Garcke, K. F. Lam, E. Sitka, V. Styles,
A Cahn--Hilliard--Darcy model for tumour growth with chemotaxis and active transport,
\textit{Math. Models Methods Appl. Sci.} \textbf{26} (2016), 1095--1148.

\bibitem{GioGrWu}
A. Giorgini, M. Grasselli, H. Wu,
The Cahn--Hilliard--Hele--Shaw system with singular potential,
{\it Ann. Inst. H. Poincar\'e Anal. Non Lin\'eaire}, {\bf 35} (2018), 1079--1118.

\pcol{\bibitem{G-B}
R. Granero-Belinch\'on, 
Global solutions for a hyperbolic-parabolic system of chemotaxis,
{\it J. Math. Anal. Appl.}  {\bf 449} (2017), 872--883.}

\bibitem{HDPZO}
A. Hawkins-Daarud, S. Prudhomme, K.~G. van der Zee, J.~T. Oden,
Bayesian calibration, validation, and uncertainty quantification of diffuse 
interface models of tumor growth,
{\it J. Math. Biol.} {\bf 67} (2013), 1457--1485.

\bibitem{HZO12}
A. Hawkins-Daarud, K. G. van der Zee, J. T. Oden,
Numerical simulation of a thermodynamically consistent four-species tumor growth model,
\textit{Int. J. Numer. Meth. Biomed. Engrg.} \textbf{28} (2012), 3--24.

\pcol{\bibitem{INK}
R. W. Ibrahim, H. K. Nashine, N. Kamaruddin, 
Hybrid time-space dynamical systems of growth bacteria with applications in segmentation,
{\it Math. Biosci.}  {\bf 292}  (2017), 10--17.}

\bibitem{JWZ}
J. Jiang, H. Wu, S. Zheng,
Well-posedness and long-time behavior of a non-autonomous
Cahn--Hilliard--Darcy system with mass source modeling tumor growth,
{\it J. \pier{Differential} Equations} {\bf 259} (2015), 3032--3077.

\pcol{\bibitem{JJ}
H. Joshi, B. K. Jha, 
Fractionally delineate the neuroprotective function of 
calbindin-D28k in Parkinson's disease,
{\it Int. J. Biomath.}  {\bf 11}  (2018),  1850103, 19 pp.}

\pcol{\bibitem{KU}		
K. H. Karlsen, S. Ulusoy, 
On a hyperbolic Keller--Segel system with degenerate nonlinear fractional diffusion,
{\it Netw. Heterog. Media}  {\bf 11} (2016), 181--201.}

\bibitem{Lions}
J.-L. Lions, ``Quelques M\'ethodes de R\'esolution des Probl\`emes aux Limites non Lin\'eaires'', 
\pier{Dunod}, Gauthier-Villars, Paris, 1969. 

\bibitem{LTZ}
J. S. Lowengrub, E. S. Titi, K. Zhao,
Analysis of a mixture model of tumor growth,
{\it European J. Appl. Math.} \textbf{24} (2013), 691--734.


\pcol{\bibitem{MV}
A. Massaccesi, E. Valdinoci, 
Is a nonlocal diffusion strategy convenient for biological populations in competition?,
{\it J. Math. Biol.}  {\bf 74}  (2017), 113--147.}

\bibitem{MRS}
A. Miranville, E. Rocca, G. Schimperna,
On the long time behavior of a tumor growth model,
{\it J. Differential Equations\/} {\bf 267} (2019) 2616--2642.

\bibitem{OHP}
J.~T. Oden, A. Hawkins, S. Prudhomme,
General diffuse-interface theories and an approach to predictive tumor growth modeling,
{\it Math. Models Methods Appl. Sci.} {\bf 20} (2010), 477--517. 

\bibitem{S}
A. Signori,
Optimal distributed control of an extended model of tumor 
growth with logarithmic potential,
{\it Appl. Math. Optim.} (2018), Online First 30 October 2018, 
https://doi.org/10.1007/s00245-018-9538-1.

\bibitem{S_DQ}
A. Signori,
Optimality conditions for an extended tumor growth model with 
double obstacle potential via deep quench approach, 
preprint arXiv:1811.08626 [math.AP] (2018), pp. 1--25.

\bibitem{S_b}
A. Signori,
Optimal treatment for a phase field system of Cahn--Hilliard 
type modeling tumor growth by asymptotic scheme, 
preprint arXiv:1902.01079 [math.AP] (2019), pp. 1--28.

\bibitem{S_a}
A. Signori,
Vanishing parameter for an optimal control problem modeling tumor growth,
preprint arXiv:1903.04930 [math.AP] (2019), pp. 1--22.

\pier{\bibitem{Simon}
J. Simon,
Compact sets in the space $L^p(0,T; B)$,
{\it Ann. Mat. Pura Appl.~(4)\/} 
{\bf 146} (1987), 65--96.}

\pcol{\bibitem{SAJM}
A. Sohail, S. Arshad, S. Javed, K. Maqbool, 
Numerical analysis of fractional-order tumor model,
{\it Int. J. Biomath.}  {\bf 8} (2015), 1550069, 13 pp.}

\bibitem{SW}
J. Sprekels, H. Wu, Optimal distributed control of a 
Cahn--Hilliard--Darcy system with mass sources,
{\em Appl. Math. Optim.} (2019), Online First 24 January 2019,
https://doi.org/10.1007/s00245-019-09555-4.

\pcol{\bibitem{SA}
N. H. Sweilam, S. M. Al-Mekhlafi,   
Optimal control for a nonlinear mathematical model of tumor under 
immune suppression: a numerical approach,
{\it Optimal Control Appl. Methods\/} {\bf 39} (2018), 1581--1596.}

\bibitem{WW2012}
X.-M. Wang, H. Wu,
Long-time behavior for the Hele--Shaw--Cahn--Hilliard system,
\textit{Asymptot. Anal.} \textbf{78} (2012), 217--245.

\bibitem{WZ2013}
X.-M. Wang, Z.-F. Zhang,
Well-posedness of the Hele--Shaw--Cahn--Hilliard system,
\textit{Ann. Inst. H. Poincar\'e Anal. Non Lin\'eaire} \textbf{30} (2013), 367--384.

\bibitem{Wise2011}
S. M. Wise, J. S. Lowengrub, V. Cristini,
An adaptive multigrid algorithm for simulating solid tumor growth using mixture models,
\textit{Math. Comput. Modelling} \textbf{53} (2011), 1--20.

\bibitem{Lowen08}
S. M. Wise, J. S. Lowengrub, H. B. Frieboes, V. Cristini,
Three-dimensional multispecies nonlinear tumor growth - I: model and numerical method,
{\it J. Theoret. Biol.} {\bf 253} (2008), 524--543.

\bibitem{WZZ}
X. Wu, G.~J. van Zwieten, K.~G. van der Zee, Stabilized second-order splitting
schemes for Cahn--Hilliard~models with applications to 
diffuse-interface tumor-growth models, 
{\it Int. J. Numer. Methods Biomed. Eng.} {\bf 30} (2014), 180--203.

\pcol{\bibitem{ZSMD}
Y. Zhou, L. Shangerganesh, J. Manimaran, A. Debbouche, 
A class of time-fractional reaction-diffusion equation with nonlocal boundary condition,
{\it Math. Methods Appl. Sci.}  {\bf 41} (2018), 2987--2999.}

\End{thebibliography}

}

\End{document}
